\documentclass[11pt]{amsart}
\usepackage{amsmath,amssymb,amsthm,amscd,verbatim,graphics}
\numberwithin{equation}{section}

\newtheorem{thm}[equation]{Theorem}
\newtheorem{lem}[equation]{Lemma}

\newtheorem{prop}[equation]{Proposition}
\newtheorem{cor}[equation]{Corollary}
\newtheorem*{thm*}{Theorem}

\theoremstyle{definition}
\newtheorem*{defn}{Definition}
\newtheorem{example}[equation]{Example}

\theoremstyle{remark}
\newtheorem{remark}[equation]{Remark}

  \newcommand{\ZZ}{\mathbb{Z}}
  \newcommand{\RR}{\mathbb{R}}
    \newcommand{\PP}{\mathbb{P}}
        
    \newcommand{\CC}{\mathbb{C}}
  \renewcommand{\H}{\mathbf{H}}

\newcommand{\Aut}{\operatorname{Aut}}
\newcommand{\Div}{\operatorname{Div}}
\newcommand{\ddiv}{\operatorname{div}}
\newcommand{\Jac}{\operatorname{Jac}}
\newcommand{\Prin}{\operatorname{Prin}}

\newcommand{\GL}{\operatorname{GL}}

\renewcommand{\setminus}{\smallsetminus}

\newcommand{\outdeg}{\operatorname{outdeg}}

\begin{document}

\title {Harmonic morphisms and hyperelliptic graphs}

\date{July 6, 2007}

\author{Matthew Baker}
\address{School of Mathematics\\
Georgia Institute of Technology\\
Atlanta, Georgia 30332-0160\\
USA}
\email{mbaker@math.gatech.edu \\ snorine@math.gatech.edu}

\author{Serguei Norine}

\begin{abstract}
We study harmonic morphisms of graphs as a natural discrete analogue of holomorphic maps between Riemann surfaces.
We formulate a graph-theoretic analogue of the classical Riemann-Hurwitz formula, study the functorial maps on Jacobians and harmonic 1-forms
induced by a harmonic morphism, and present a discrete analogue of the canonical map from a Riemann surface to projective space.
We also discuss several equivalent formulations of the notion of a hyperelliptic graph, all motivated by the classical theory of Riemann surfaces.
As an application of our results, we show that for a 2-edge-connected graph $G$ which is not a cycle, there is at most one involution $\iota$ on $G$ for
which the quotient $G/\iota$ is a tree.  We also show that the number of spanning trees in a graph $G$ is even if and only if $G$ admits a non-constant 
harmonic morphism to the graph $B_2$ consisting of 2 vertices connected by 2 edges.  
Finally, we use the Riemann-Hurwitz formula and our results on hyperelliptic graphs to 
classify all hyperelliptic graphs having no Weierstrass points.
\end{abstract}

\maketitle

\section{Introduction}

\subsection{Notation and terminology}
\label{NotationSection}

Throughout this paper, a {\em Riemann surface} will mean a compact, connected one-dimensional complex manifold, and
(unless otherwise specified) a {\em graph} will mean a finite, connected multigraph without loop edges.
A graph with no multiple edges is called {\em simple}.
We will denote by $V(G)$ and $E(G)$, respectively, the set of vertices and edges of $G$.
For a vertex $x \in V(G)$ and an edge $e \in E(G)$, we write $x \in
e$ if $e$ is incident to $x$.

We denote by $g(G) := |E(G)| - |V(G)| + 1$ 
the {\em genus} of $G$; this is the dimension of the cycle space of 
$G$.
(In graph theory, the term ``genus'' is traditionally used for a different concept, namely, the smallest genus of any surface in which the graph can be embedded, and the integer $g = g(G)$ is called the ``cyclomatic number'' of $G$.  We call $g$ the genus of $G$ in order to highlight the analogy with Riemann surfaces.)

For $k \geq 2$, a graph $G$ is called {\em $k$-edge-connected} if
$G-W$ is connected for every set $W$ of at most $k-1$ edges of $G$.
(By convention, we consider the trivial graph having one vertex and no edges to be
$k$-edge-connected for all $k$.)
Alternatively, define a {\em cut} to be the set of all edges connecting a vertex in $V_1$ to a vertex in $V_2$
for some partition
of $V(G)$ into disjoint subsets $V_1$ and $V_2$.
Then $G$ is $k$-edge-connected if and only if every non-empty cut has size at least $k$.

A {\em bridge} is an edge of $G$ 
whose deletion increases the number of connected components of $G$.
A (connected) graph is $2$-edge-connected if and only if it contains no bridge.

Finally, if $A \subseteq V(G)$, we denote by $\chi_A : V(G) \to \{ 0,1 \}$ the characteristic function
of $A$, and for $x \in A$ we let $\outdeg_A(x)$ denote the number of edges $e=xy$ in $E(G)$ with $y \not\in A$.

\subsection{Motivation and discussion of main results}

In \cite{BakerNorine}, the authors investigated some new analogies between graphs and Riemann surfaces,
formulating the notion of a {\em linear system} on a graph and proving a graph-theoretic analogue of the classical
Riemann-Roch theorem.
The theory of linear systems on graphs has applications to understanding the {\em Jacobian of a finite graph}, 
a group which is analogous to the Jacobian of a Riemann surface, and which 
has appeared in many different guises throughout the literature (e.g., as the ``Picard group'' in \cite{BDN}, the 
``critical group'' in \cite{BiggsAPTG}, the ``sandpile group'' in \cite{DharSandpile}, and the
``group of components'' in \cite{Lorenzini2}).

\medskip

The present paper can be viewed as a natural sequel to \cite{BakerNorine}.
In classical algebraic geometry,
one is usually interested not just in Riemann surfaces themselves, 
but also in the holomorphic maps between them.  
Thus, we are naturally led to ask: what is the ``correct'' graph-theoretic analogue of a holomorphic map
between Riemann surfaces?  In other words, is there a category consisting of graphs, together with certain maps between them,
which closely mirrors the category of Riemann surfaces, together with the holomorphic maps between them?
In this paper, we hope to convince the reader 
that the notion of a {\em harmonic morphism of graphs}, introduced by Urakawa in \cite{Urakawa1},
has essentially all of the desired features.

\medskip

Actually, since we want to allow graphs with multiple edges, we need to slightly modify the definition of a
harmonic morphism from \cite{Urakawa1}, since Urakawa assumes that all of his graphs are simple.
Recall that a holomorphic map $\phi : X \to X'$ between Riemann surfaces is one which locally 
pulls back holomorphic functions on $X'$ to holomorphic functions on $X$.
Although the notion of a holomorphic function does not make sense in the context of graphs, there
is a natural notion of a {\em harmonic function}
(see (\ref{e:harmonicmapdefn}) below for the definition).
Urakawa defines a harmonic morphism $\phi : V(G) \to V(G')$ between simple graphs $G,G'$ 
to be a function which locally pulls back harmonic functions on $G'$ to harmonic functions on $G$, 
i.e., a function such that for every $x \in V(G)$ and every function $f : V(G') \to \RR$
which is harmonic at $\phi(x)$, the composition $f \circ \phi$ is harmonic at $x$.
What makes this a workable and useful notion is Theorem 2.5 from \cite{Urakawa1}, which asserts
that a function $\phi : V(G) \to V(G')$ between simple graphs is a harmonic morphism if and only if $\phi$ is {\em horizontally conformal},
meaning that:
\begin{itemize}
\item[(HC1)] For all adjacent vertices $x,y \in V(G)$, we have either $\phi(x) = \phi(y)$ or $\phi(x)$ is adjacent to $\phi(y)$, and
\item[(HC2)] For all $x \in V(G)$, the quantity 
\[
|\{ y \in V(G) \; | \; \phi(y) = y' \textrm{ and } y \textrm{ is adjacent to } x \}|
\]
is the same for all $y' \in V(G')$ adjacent to $x' := \phi(x)$.
\end{itemize}

However, the equivalence between harmonic morphisms and horizontally conformal maps fails 
for graphs which are not simple (c.f.~Remark~\ref{NotSimpleRemark} below).
Because of this, we take an analogue of (HC1) and (HC2) as our definition of a harmonic morphism 
between multigraphs (see Definition~\ref{HarmonicMorphismDefinition} for a precise definition).  
A harmonic morphism in this sense does indeed pull back harmonic functions to harmonic functions,
but the converse does not always hold.  

\medskip

One of the key features of defining harmonic morphisms in terms of horizontal conformality
is that given a harmonic morphism $\phi : G \to G'$, it is possible to assign a well-defined {\em multiplicity}
$m_{\phi}(x)$ to each vertex $x \in V(G)$ (analogous to the {\em ramification index} $e_{\phi}(x)$ at $x \in X$ for a non-constant holomorphic map 
$\phi : X \to X'$ between Riemann surfaces) in such a way that the sum $\deg(\phi)$ of the multiplicities at all vertices mapping to a given vertex
$x' \in V(G')$ is independent of $x'$.  
We define the {\em degree} of $\phi$ to be this number.

\medskip

Harmonic morphisms between graphs enjoy numerous properties analogous to classical properties from algebraic geometry.
For example, if $\phi : X \to X'$ is a non-constant holomorphic map of degree $\deg(\phi)$ between Riemann surfaces having genus $g$ and $g'$, respectively, 
then:

\begin{itemize}
\item[(RS1)] $\phi$ is surjective and $g \geq g'$.
\item[(RS2)] The {\em Riemann-Hurwitz} formula $2g-2 = \deg(\phi)(2g'-2) + \sum_{x \in X} \left( e_{\phi}(x) - 1 \right)$ holds.
\item[(RS3)] $\phi$ induces functorial maps $\phi_* : \Jac(X) \to \Jac(X')$ and $\phi^* : \Jac(X') \to \Jac(X)$ between the Jacobians of
$X$ and $X'$.
\item[(RS4)] $\phi$ induces functorial maps $\phi_* : \Omega^1(X) \to \Omega^1(X')$ and $\phi^* : \Omega^1(X') \to \Omega^1(X)$ between the complex
vector spaces $\Omega^1(X)$ and $\Omega^1(X')$ of holomorphic $1$-forms on $X$ and $X'$, respectively.
\item[(RS5)] If $D$ is a divisor on $X$, then $\dim |\phi_*(D)| \geq \dim |D|$, where $|D|$ denotes the complete linear system associated to $D$.  
In particular, if $X$ is hyperelliptic and $g(X') \geq 2$, then $X'$ is hyperelliptic as well.
\end{itemize}

We will prove graph-theoretic analogues of all of these classical facts.  
We will also describe some situations in which the naive analogue of certain classical facts does not hold.
For example, in algebraic geometry the map $\phi^* : \Jac(X') \to \Jac(X)$ is sometimes injective and sometimes not; 
more precisely, it is known that $\phi^*$ fails to be injective if and only if $\phi$ has a nontrivial unramified abelian subcover.
However, the analogous map $\phi^* : \Jac(G') \to \Jac(G)$ in the graph-theoretic context turns out to always be injective; 
this appears to be a rather subtle fact with some useful applications.

\medskip

As a basic testing ground for our ``dictionary'' between graphs and Riemann surfaces, 
we consider in detail the graph-theoretic analogue of
a {\em hyperelliptic} Riemann surface.  
This is particularly interesting because classically, there are many different equivalent characterizations
of what it means for a Riemann surface $X$ of genus at least $2$ to be hyperelliptic.
As just a sample, we mention the following:

\begin{itemize}
\item[(H1)] There exists a divisor $D$ of degree $2$ on $X$ for which $r(D):= \dim |D|$ is equal to $1$.
\item[(H2)] There exists an involution $\iota$ for which $X/\iota$ has genus $0$.
\item[(H3)] There is a degree $2$ holomorphic map $\phi : X \to \PP^1$.
\item[(H4)] There is an automorphism $\iota$ of $X$ for which 
$\iota^* : \Jac(X) \to \Jac(X)$ is multiplication by $-1$.
\item[(H5)] There is an automorphism $\iota$ of $X$ for which 
$\iota^* : \Omega^1(X) \to \Omega^1(X)$ is multiplication by $-1$.  
\item[(H6)] The symmetric square $S_{x_0}^{(2)} : \Div_{+}^{2}(X) \to \Jac(X)$ of the Abel-Jacobi map
(relative to some base point $x_0 \in X$) is not injective.
\item[(H7)] The canonical map $\psi_X : X \to \PP(\Omega^1(X))$ is not injective.
\end{itemize}

When any of these equivalent conditions are satisfied, there is a unique automorphism $\iota$ satisfying 
(H2), (H4), and (H5), called the {\em hyperelliptic involution}.

\medskip

For a $2$-edge-connected graph $G$ of genus at least $2$, we take the analogue of (H1) 
to be the definition of what it means for $G$ to be hyperelliptic.
(This definition was already introduced in \cite{BakerSpecialization}.)
We then prove that the graph-theoretic analogues of conditions (H1)-(H5) above are all equivalent to one another,
and that the hyperelliptic involution $\iota$ on a graph satisfying any of these conditions is unique.
However, in the graph-theoretic context it turns out that $\textrm{(H1)} \Rightarrow \textrm{(H6)}  
\Leftrightarrow \textrm{(H7)}$, so that hyperelliptic graphs satisfy the analogues of 
conditions (H6) and (H7), but there are non-hyperelliptic $2$-edge-connected graphs $G$ 
of genus at least $2$ which also satisfy these conditions.
In fact, we will see that the graph-theoretic analogues of conditions (H6) and (H7) are equivalent to the condition 
that $G$ is not $3$-edge-connected.

\medskip

As an application of our results, and to illustrate another difference with the theory of Riemann surfaces, 
we conclude our paper with a discussion of Weierstrass points on hyperelliptic graphs.
(The notion of a Weierstrass points on graphs was introduced in \cite{BakerSpecialization}; see 
\S\ref{HyperellipticSection} for a definition.)
Classically, a hyperelliptic Riemann surface of genus $g\geq 2$ possesses exactly $2g+2$ Weierstrass points,
namely, the fixed points of the hyperelliptic involution, and every Riemann surface of genus at least $2$
has Weierstrass points.  The situation for graphs is less orderly,
as there are infinite families of graphs having no Weierstrass points at all.  
Using our rather precise knowledge about the structure of hyperelliptic graphs,
we give a classification of all hyperelliptic graphs having no Weierstrass points.
We leave as an open problem whether or not there exist further (non-hyperelliptic) examples
of Weierstrass-pointless graphs.

\medskip

Occasionally, our foundational results on harmonic morphisms and hyperelliptic graphs lead to 
applications to more traditional-sounding graph-theoretic questions.
For example, as a consequence of our study of hyperelliptic graphs, we will show that for a 
$2$-edge-connected graph $G$ of genus at least $2$, there is at most one involution $\iota$ on $G$
whose quotient is a tree.  We also show that the number $\kappa_G$ of spanning trees in a graph $G$ is even 
if and only if $G$ admits a non-constant degree $2$ harmonic morphism to the graph $B_2$ consisting
of $2$ vertices connected by $2$ edges.

\medskip

Although in this paper we view our graph-theoretic results as ``analogous'' to classical results from algebraic geometry,
there is in fact a closer connection between the two worlds than one might at first imagine.  One such connection
arises from the specialization of divisors on an arithmetic surface, and is explored in 
\cite{BakerSpecialization}.  
We expect that the ideas in the present paper will help spur further interactions between graph theory, on the one hand, and 
arithmetic, algebraic, and tropical geometry on the other.

\medskip

It would be interesting to prove analogues of the results in the present paper for metric graphs, and more generally for tropical curves,
but we have not attempted to do so here.  It would also be interesting to 
generalize some of our results to higher dimensions.  At least in the context of Riemannian polyhedra (which are
higher-dimensional analogues of metric graphs), there is already a rich literature concerning
the notion of a harmonic morphism (see, e.g., \cite{FugledeBook}).  However, it appears that the questions 
being addressed in \cite{FugledeBook} and the references therein are of a somewhat different flavor
than the ones which we study here.

\medskip

We have endeavored to make this paper as self-contained as possible.  Therefore, we summarize in \S\ref{GraphDivisorSection}
below all of the facts from \cite{BakerNorine} which we will be using.  We have also rewritten certain proofs
from \cite{Urakawa1}, because our notation differs somewhat from Urakawa's, and because we work
in the somewhat more general setting of multigraphs.  
A good reference for many of the facts about Riemann surfaces which we will be discussing in this paper is \cite{Miranda},
and a basic graph theory reference is \cite{Bollobas}.

\subsection{Background material from \cite{BakerNorine}}
\label{GraphDivisorSection}

Following \cite{BakerNorine}, we denote by $\Div(G)$ 
the free abelian group on $V(G)$.  
We refer to elements of $\Div(G)$ as {\em divisors} on $G$.  
We can write each element $D \in \Div(G)$ uniquely as
$D = \sum_{x \in V(G)} D(x) (x)$ with $D(x) \in \ZZ$.
We say that $D$ is {\em effective}, and write $D \geq 0$, 
if $D(x) \geq 0$ for all $x \in V(G)$.
For $D \in \Div(G)$, we define the {\em degree} of $D$ by the formula
$\deg(D) = \sum_{x \in V(G)} D(x)$.
We denote by $\Div^0(G)$ the subgroup of $\Div(G)$ consisting of divisors of degree zero.
Finally, we denote by $\Div_{+}^{k}(G) = \{ E \in \Div(G) \; : \; E \geq 0, \, \deg(E)=k \}$ 
the set of effective divisors of degree $k$ on $G$.

Let $C^0(G,\ZZ)$ be the group of $\ZZ$-valued functions on $V(G)$.
For $f \in C^0(G,\ZZ)$, we define the {\em divisor} of $f$ by the formula
\[
\ddiv(f) = \sum_{x \in V(G)} \sum_{e = xy \in E(G)} 
\left( f(x) - f(y) \right)(x).
\]

The divisor of $f$ can be naturally identified with the graph-theoretic Laplacian of $f$.
Divisors of the form $\ddiv(f)$ for some $f \in C^0(G,\ZZ)$ are called {\em principal}; we
denote by $\Prin(G)$ the group of principal divisors on $G$.
It is easy to see that every principal divisor has degree zero, so that
$\Prin(G)$ is a subgroup of $\Div^0(G)$.

The Jacobian of $G$, denoted $\Jac(G)$, is defined to be the quotient group
\[
\Jac(G) = \Div^0(G) / \Prin(G).
\]

One can show using Kirchhoff's Matrix-Tree Theorem (c.f.~\cite[\S{14}]{BiggsAPTG}) that $\Jac(G)$ is a finite abelian group of 
order $\kappa_G$, where $\kappa_G$ is the number of spanning trees in $G$.

\medskip

We define an equivalence relation $\sim_G$ on $\Div(G)$ by writing
$D \sim_G D'$ if and only if $D - D' \in \Prin(G)$,
and set
\[
|D| = \{ E \in \Div(G) \; : \; E \geq 0 \textrm{ and } E \sim_G D \}.
\]
We refer to $|D|$ as the {\em (complete) linear system} associated
to $D$, and when $D \sim D'$ we call the divisors $D$ and $D'$ {\em linearly equivalent}.
We will usually just write $D \sim D'$, rather than $D \sim_G D'$, when the graph $G$ is
understood.

For later use, we note the following simple fact about the linear equivalence relation on $G$:

\begin{lem}
\label{l:allequiv}
We have $(x) \sim (y)$ for all $x,y \in V(G)$ if and only if $G$ is a tree.
\end{lem}

\begin{proof}
This follows from the fact that $|\Jac(G)| = \kappa_G$, together with the observation
that the group $\Div^0(G)$ is generated by the divisors of the form $(x)-(y)$ with $x,y \in V(G)$.
\end{proof}

\medskip

Given a divisor $D$ on $G$, define
$r(D) = -1$ if $|D| = \emptyset$, and otherwise set
\[
r(D) = \max \{ k \in \ZZ \; : \; |D - E| \neq \emptyset \; \forall \; E \in \Div_+^{k}(G) \}.
\]
Note that $r(D)$ depends only on the linear equivalence class of $D$,
and therefore is an invariant of the linear system $|D|$.
When we wish to emphasize the underlying graph $G$, we will
sometimes write $r_G(D)$ instead of $r(D)$.

For later use, we recall from \cite[Lemma~2.1]{BakerNorine} the following simple lemma:

\begin{lem}
\label{SubAdditivityLemma}
For all $D,D' \in \Div(G)$ such that
$r(D),r(D') \geq 0$, we have $r(D + D') \geq r(D) + r(D')$.
\end{lem}

\medskip

We define the {\em canonical divisor} on $G$ to be 
\[
K_G = \sum_{x \in V(G)} (\deg(x) - 2)(x).
\]
We have $\deg(K_G) = 2g-2$, where $g = |E(G)| - |V(G)| + 1$ 
is the {\em genus} of $G$ (or, in more traditional language, {\em cyclomatic number} of $G$).

\medskip

The following result is proved in \cite[Theorem 1.12]{BakerNorine}:

\begin{thm}[Riemann-Roch for graphs]
\label{t:GraphRiemannRoch}
Let $D$ be a divisor on a graph $G$.  Then
\[
r(D) - r(K_G - D) = \deg(D) + 1 - g.
\]
\end{thm}

As a consequence of Lemma~\ref{SubAdditivityLemma} and Theorem~\ref{t:GraphRiemannRoch},
we have the following graph-theoretic analogue of a classical result known as Clifford's theorem
(see \cite[Corollary~3.5]{BakerNorine} for a proof):

\begin{cor}[Clifford's Theorem for graphs]
\label{CliffordCor}
Let $D$ be a divisor on a graph $G$ for which $|D| \neq \emptyset$ and $|K_G-D| \neq \emptyset$.  Then
\[
r(D) \leq \frac{1}{2} \deg(D) \ .
\]
\end{cor}

The next result (Theorem~3.3 from \cite{BakerNorine}) is very useful for computing
$r(D)$ in specific examples, and also plays an important role in the proof of Theorem~\ref{t:GraphRiemannRoch}.  
For each linear ordering $<$ on $V(G)$,
we define a corresponding divisor $\nu \in \Div(G)$ of degree $g-1$
by the formula
\[
\nu = \sum_{x \in V(G)}(|\{e = xy \in E(G) \; : \; y < x \}| -
1)(x).
\]

\begin{thm}
\label{t:OrderDivisor}
For every $D \in \Div(G)$, exactly one of the following holds:
\begin{enumerate}
\item $r(D) \ge 0$; or
\item $r(\nu - D) \ge 0$ for some divisor $\nu$ associated to a linear ordering $<$ of $V(G)$.
\end{enumerate}
\end{thm}

Finally, we recall some facts from \cite{BakerNorine} and \cite{BDN} about the graph-theoretic analogue of the 
Abel-Jacobi map from a Riemann surface to its Jacobian.

For a fixed base point $x_0 \in V(G)$, we define the {\em Abel-Jacobi map} $S_{x_0} : G \to \Jac(G)$
by the formula
\begin{equation}
\label{eq:AbelJacobiMap}
S_{x_0}(x) = [(x) - (x_0)] \ .
\end{equation}

The map $S_{x_0}$ can be characterized by the following universal property (see \S{3} of \cite{BDN}).
A map $\varphi : G \to A$ from $V(G)$ to an abelian group $A$ is called {\em harmonic} if
for each $x \in V(G)$, we have
\begin{equation}
\label{e:harmonicmapdefn}
\deg(x) \cdot \varphi(x) = \sum_{e=xy \in E(G)} \varphi(y) \ .
\end{equation}
Then $S_{x_0}$ is universal among all harmonic maps from $G$ to abelian groups sending $x_0$ to $0$,
in the following precise sense:

\begin{lem}
\label{l:AbelJacobiUnivProperty}
If $\varphi : G \to A$ is any map sending $x_0 \in V(G)$ to $0$, then there is a unique group homomorphism
$\psi : \Jac(G) \to A$ such that $\varphi = \psi \circ S_{x_0}$.
\end{lem}

We also define, for each integer $k \geq 1$, a map $S_{x_0}^{(k)} : \Div_+^k(G) \to \Jac(G)$
by
\[
S_{x_0}^{(k)}((x_1) + \cdots + (x_k)) = S_{x_0}(x_1) + S_{x_0}(x_2) + \cdots + S_{x_0}(x_k) \ .
\]

The following result is proved in \cite[Theorem 1.8]{BakerNorine}:

\begin{thm}
\label{t:GraphAbelJacobi}
The map $S_{x_0}^{(k)}$ is 
injective if and only if $G$ is $(k+1)$-edge-connected.
\end{thm}

\section{Harmonic morphisms}
\label{HarmonicMorphismSection}

\subsection{Definition and basic properties of harmonic morphisms}

Harmonic morphisms between simple graphs were defined and studied
in~\cite{Urakawa1}. Here, we reproduce some definitions
from~\cite{Urakawa1}, but with several variations due to the fact that we allow
our graphs to have multiple edges.

\medskip

Let $G, G'$ be graphs. A function $\phi: V(G) \cup E(G) \rightarrow
V(G') \cup E(G')$ is said to be a \emph{morphism} from $G$ to $G'$
if $\phi(V(G)) \subseteq V(G')$, and for every $x \in V(G)$ and $e \in
E(G)$ such that $x \in e$, either $\phi(e)\in E(G')$ and $\phi(x) \in
\phi(e)$, or $\phi(e)=\phi(x)$. We write $\phi: G \rightarrow G'$ for brevity.
If $\phi(E(G)) \subseteq E(G')$ then we say that $\phi$ is a
\emph{homomorphism}.  A bijective homomorphism is called an
\emph{isomorphism}, and an isomorphism $\phi: G \rightarrow G$ 
is called an \emph{automorphism}.

\medskip

We now come to the key definition in this paper.

\begin{defn}
\label{HarmonicMorphismDefinition}
A morphism $\phi: G \rightarrow G'$ is said to be \emph{harmonic}
(or \emph{horizontally conformal}) if for
all $x \in V(G),y\in V(G')$ such that $y = \phi(x)$, the quantity
$|\{e \in E(G) |x \in e,\; \phi(e) =e' \}|$ is the same for all edges $e' \in E(G')$
such that $y \in e'$.
\end{defn}

\begin{remark}
One can check directly from the definition that the composition of two harmonic morphisms is again harmonic.
Therefore the set of all graphs, together with the harmonic morphisms between them, forms a category.
\end{remark}

Let $\phi: G \rightarrow G'$ be a morphism and let $x \in
V(G)$.
Define the {\em vertical multiplicity} of $\phi$ at $x$ by
\[
v_{\phi}(x)= |\{e \in E(G) \; | \phi(e) =\phi(x)\}|.
\]
This is simply the number of {\em vertical edges} incident to $x$, where an edge $e$ is
called {\em vertical} if $\phi(e) \in V(G)$ (and is called {\em horizontal} otherwise).

\medskip

If $\phi$ is harmonic and $|V(G')| > 1$, we define the {\em horizontal
multiplicity} 
of $\phi$ at $x$ by
\[
m_{\phi}(x) = |\{e \in E(G)  |x \in e, \; \phi(e) =e'  \}|
\]
for any edge $e' \in E(G)$ such that $\phi(x) \in e'$.
By the definition of a harmonic morphism,
$m_{\phi}(x)$ is independent of the choice of $e'$.
When $|V(G')| = 1$, we define $m_{\phi}(x)$ to be $0$ for all $x \in V(G)$.

\medskip

If $\deg(x)$ denotes the degree of a vertex $x$, we have the following basic formula
relating the horizontal and vertical multiplicities:

\begin{equation}\label{e:DegreeIdentity}
\deg(x) = \deg(\phi(x))m_{\phi}(x)+v_{\phi}(x).
\end{equation}

\medskip

We say that a harmonic morphism $\phi: G \rightarrow G'$ is
\emph{non-degenerate} if $m_{\phi}(x) \geq 1$ for every $x \in
V(G)$.  (The motivation for this definition comes from
Theorem~\ref{t:Hyperelliptic} below.)

\medskip

If $|V(G')| > 1$, we define the {\em degree} of a harmonic morphism $\phi: G \rightarrow G'$ by the formula
\begin{equation}
\label{e:DegreeDef}
\deg(\phi) := |\{ e \in E(G) \; | \; \phi(e) = e'\}|
\end{equation}
for any edge $e' \in E(G')$.  (When $|V(G')| = 1$, we define $\deg(\phi)$ to be $0$.)
By the following lemma (c.f.~\cite[Lemma 2.12]{Urakawa1}),
the right-hand side of (\ref{e:DegreeDef}) does not depend on the choice of $e'$
(and therefore $\deg(\phi)$ is well-defined):

\begin{lem}
\label{lem:edgeindependence} The quantity $|\{e \in E(G) \; | \;
\phi(e) = e' \}|$ is independent of the choice of $e' \in E(G')$.
\end{lem}

\begin{proof}
Let $y \in V(G')$, and suppose there are two edges $e',e'' \in E(G')$ incident to $y$.
Since $\phi$ is horizontally conformal, for each $x \in V(G)$ with $\phi(x) = y$ we have
\[
| \{ e \in E(G) \; | \; x \in e, \, \phi(e) = e' \}|
= | \{ \tilde{e} \in E(G) \; | \; x \in \tilde{e}, \, \phi(\tilde{e}) = e'' \}|.
\]
Therefore
\begin{equation}
\label{e:edgecalc}
\begin{aligned}
|\{ e \in E(G) \; | \; \phi(e) = e'\} | &=
\sum_{x \in \phi^{-1}(y)} | \{ e \in E(G) \; | \; x \in e, \, \phi(e) = e' \}| \\
&= \sum_{x \in \phi^{-1}(y)} | \{ \tilde{e} \in E(G) \; | \; x \in \tilde{e}, \, \phi(\tilde{e}) = e'' \} | \\
&= | \{ \tilde{e} \in E(G) \; | \; \phi(\tilde{e}) = e'' \} |. \\
\end{aligned}
\end{equation}

Now suppose $e',e''$ are arbitrary edges of $G'$.
Since $G$ is connected, the result follows by applying (\ref{e:edgecalc}) to each pair of consecutive
edges in any path connecting $e'$ and $e''$.
\end{proof}

\medskip

According to the next result, the degree of a harmonic morphism $\phi : G \to G'$ is just the number
of preimages under $\phi$ of any vertex of $G'$, counting multiplicities:

\begin{lem}
\label{lem:vertexindependence}
For any vertex $y \in G'$, we have
\[
\deg(\phi) = \sum_{\substack{x \in V(G) \\ \phi(x) = y}} m_{\phi}(x).
\]
\end{lem}

\begin{proof}
Choose an edge $e' \in E(G')$ with $y \in e'$.  Then
\[
\begin{aligned}
\sum_{x \in \phi^{-1}(y)} m_{\phi}(x) &= \sum_{x \in \phi^{-1}(y)} \sum_{e \in \phi^{-1}(e'), \, x \in e} 1 \\
&= | \phi^{-1}(e') | = \deg(\phi). \\
\end{aligned}
\]
\end{proof}

As with morphisms of Riemann surfaces in algebraic geometry, a harmonic morphism of graphs must be
either constant or surjective.   More generally, one has the following:

\begin{lem}
\label{lem:constantorsurj}
Let $\phi : G \rightarrow G'$ be a harmonic morphism with $|V(G')|>1$.  
Then $\deg(\phi) = 0$ if and only if $\phi$ is constant, 
and $\deg(\phi) > 0$ if and only if $\phi$ is surjective.
\end{lem}

\begin{proof}
If $\phi$ is constant, then clearly $\deg(\phi) = 0$.  Moreover, it follows from Lemmas~\ref{lem:edgeindependence}
and \ref{lem:vertexindependence} that $\phi$ is surjective if and only if $\deg(\phi) > 0$.  So it remains only to
show that if $\deg(\phi) = 0$, then $\phi$ is constant.  For this, suppose we have $\phi(x) = y$.  Since $m_{\phi}(x) = 0$,
it follows that $\phi(e) = y$ for every edge $e$ with $x \in e$.  Thus $\phi(x') = y$ for every neighbor $x'$ of $x$.
As $G$ is connected, it follows that every vertex and every edge of $G$ is mapped under $\phi$ to $y$.
\end{proof}

\subsection{Harmonic morphisms and harmonic maps to abelian groups}

Recall that given a graph $G$ and an abelian group $A$, a function $f:V(G)
\rightarrow A$  is said to be \emph{harmonic} at $x\in V(G)$ if
$$\sum_{e = xy \in E(G)} (f(x)-f(y))=0.$$ A morphism $\phi: G
\rightarrow G'$ is said to be \emph{$A$-harmonic} if for any $y =
\phi(x)$ and any function $f: V(G') \rightarrow A$ harmonic at $y$,
the function $f \circ \phi$ is harmonic at $x$.

\begin{prop}
\label{p:conformity} Let $G$ and $G'$ be graphs, and let $\phi:
G \rightarrow G'$ be a harmonic morphism.  Then $\phi$ is
$A$-harmonic for every abelian group $A$.
\end{prop}

\begin{proof}(c.f.~Lemma~2.11 of \cite{Urakawa1})
Let $x \in V(G)$, $y \in V(G')$ be such that $y = \phi(x)$, and
let $f: V(G') \rightarrow A$ be harmonic at $y$, i. e.
$$\sum_{e = zy \in E(G')}
f(z)= \deg(y)f(y).$$ 
Then we have
\begin{align*}
\sum_{e = zx \in E(G)} f(\phi(z)) &= \sum_{\substack{e = zx \in E(G) \\ \phi(e)=y}}f(\phi(z)) 
+ \sum_{e' = z'y \in E(G')} \left( \sum_{\substack{e = zx \in E(G) \\ \phi(e)=e'}}f(\phi(z)) \right) \\
&= v_{\phi}(x)f(y) + \sum_{e' = z'y \in E(G')} m_{\phi}(x) f(z') \\ 
&= v_{\phi}(x)f(y) + m_{\phi}(x)\deg(y)f(y) \\ 
&= (v_{\phi}(x) + m_{\phi}(x)\deg(\phi(x)))f(\phi(x)) \\
&=\deg(x)f(\phi(x)) \qquad\qquad\qquad\qquad\qquad\qquad\text{(by (\ref{e:DegreeIdentity})),}
\end{align*}
as desired.
\end{proof}

If $G'$ is a simple graph (i.e., without multiple edges), then the converse of Proposition~\ref{p:conformity}
also holds:

\begin{prop}
\label{p:conformityconverse}
If $G'$ is a simple graph, then for a morphism $\phi: G \rightarrow G'$, the following are
equivalent:
\begin{enumerate}
  \item\label{c:Harmonic1} $\phi$ is harmonic (i.e., horizontally conformal).
  \item\label{c:Harmonic2} $\phi$ is $A$-harmonic for every abelian group $A$.
  \item\label{c:Harmonic3} $\phi$ is $\RR$-harmonic.
\end{enumerate}
\end{prop}

\begin{proof}
It follows from Proposition~\ref{p:conformity} that (\ref{c:Harmonic1})
implies (\ref{c:Harmonic2}), and it is immediate that
(\ref{c:Harmonic2}) implies (\ref{c:Harmonic3}). It remains to show
that (\ref{c:Harmonic3}) implies (\ref{c:Harmonic1}), which we do following
\cite[Lemma 2.7]{Urakawa1}.

For a vertex $x \in V(G)$ and an edge $e' \in E(G')$ such that
$\phi(x) \in e'$, let $$k_x(e')=|\{e \in E(G) |x \in e,\; \phi(e)
=e' \}|.$$ We need to prove that $k_x(e')$ is independent of the
choice of $e'$. Let $\phi(x)=y$, let $e'=yz$, and define a function
$f_{e'}: V(G') \rightarrow \RR$ as follows. Let
$f_{e'}(z)=1,$ let $f_{e'}(y)= 1/\deg(y)$, and let $f_{e'}(w)=0$ for
$w \in V(G') \backslash \{y, z\}$. Then $f_{e'}$ is harmonic at
$y$, so by (\ref{c:Harmonic3}), $f_{e'} \circ \phi$ is harmonic at $x$. It
follows that
\begin{align*}
\frac{\deg(x)}{\deg(y)}&= \deg(x) f_{e'}(\phi(x)) = \sum_{e=xw\in
E(G)}f_{e'}(\phi(w))\\ 
&= \sum_{\substack{e=xw \in E(G)\\ \phi(w)=y}}f_{e'}(y) + \sum_{\substack{e=xw \in E(G)\\ \phi(w)=z}}f_{e'}(z) \\ 
&= \frac{v_{\phi}(x)}{\deg(y)} + k_x(e') \qquad\qquad\qquad\qquad\qquad\text{(since $G'$ is simple).}
\end{align*}
Therefore $k_x(e') = (\deg(x) - v_{\phi}(x))/\deg(\phi(x))$ is
independent of the choice of $e'$,
as desired.
\end{proof}

\begin{remark}
\label{NotSimpleRemark}
If $G'$ is not simple, then the converse of Proposition~\ref{p:conformityconverse} may fail, as one
sees easily by taking $G$ to be the graph with $2$ vertices $x,y$ connected by a single edge $e$,
$G'$ to be the graph with $2$ vertices $x',y'$ connected by two edges $e'_1,e'_2$,
and $\phi : G \to G'$ to be the morphism which sends $x,y$ to $x',y'$, respectively, and
$e$ to $e'_1$.
\end{remark}

\medskip

\subsection{The Riemann-Hurwitz formula for graphs}
\label{RiemannHurwitzSection}

Let $\phi : G \to G'$ be a harmonic morphism.  We define the push-forward homomorphism
$\phi_* : \Div(G) \to \Div(G')$ by
\begin{equation}
\label{eq:Lowerstar}
\phi_*(D) = \sum_{x \in V(G)} D(x) (\phi(x)).
\end{equation}
Similarly, we define the pullback homomorphism $\phi^*: \Div(G') \rightarrow \Div(G)$ by
\begin{equation}
\label{eq:Upperstar}
\phi^*(D')= \sum_{y \in V(G')}
\sum_{\substack{x \in V(G) \\  \phi(x)=y}} m_{\phi}(x)D'(y) (x).
\end{equation}

\begin{lem}
\label{lem:PullbackDegree}
If $\phi : G \to G'$ is a harmonic morphism and $D' \in \Div(G')$, then
$\deg(\phi^*(D')) = \deg(\phi) \cdot \deg(D')$.
\end{lem}

\begin{proof}
This follows from Lemma~\ref{lem:vertexindependence} and the definition of $\phi^*$.
\end{proof}



\medskip

A basic fact about harmonic morphisms of graphs is that one has the following analogue of the classical Riemann-Hurwitz formula
from algebraic geometry:

\begin{thm}[Riemann-Hurwitz for graphs]
\label{t:RiemannHurwitz}
Let $G,G'$ be graphs,
and let $\phi: G \rightarrow G'$ be a harmonic morphism. Then:
\begin{enumerate}
\item The canonical divisors on $G$ and $G'$ are related by the formula
\begin{equation}\label{e:RiemannHurwitz}
K_{G}=\phi^*K_{G'} + R_{G},
\end{equation}
where $$R_{G}=2\sum_{x \in V(G)}(m_{\phi}(x)-1)(x) +\sum_{x \in V(G)} v_{\phi}(x)(x).$$
\item  If $G,G'$ have genus $g$ and $g'$, respectively, then
\begin{equation}
\label{e:RH2}
2g-2=\deg(\phi)(2g'-2)+\sum_{x \in V(G)}\left( 2(m_{\phi}(x)-1)+v_{\phi}(x) \right).
\end{equation}
\item If $\phi$ is non-constant, then $2g-2 \geq \deg(\phi)(2g'-2)$ and $g \geq g'$.
\end{enumerate}
\end{thm}

\begin{proof}
By definition, we have $(\phi^*K_{G'})(x)=m_{\phi}(x)(\deg(\phi(x))-2)$.
On the other hand, by (\ref{e:DegreeIdentity}) we have
\begin{align*}
K_{G}(x)&=\deg(x) -2= \deg(\phi(x))m_{\phi}(x)+v_{\phi}(x) -2 \\ &=(\phi^*K_{G'})(x) +
2m_{\phi}(x)+v_{\phi}(x) -2 = (\phi^*K_{G'}+R_{G})(x)
\end{align*}
for every $x \in V(G)$, which proves (1).
Part (2) follows immediately from Lemma~\ref{lem:PullbackDegree}
upon computing the degrees of the divisors on both sides of (\ref{e:RiemannHurwitz}).
In order to verify (3), we claim that if $\phi$ is non-constant then $\deg(R_G) \geq 0$.
This is clear if $G$ has no {\em vertical leaves} (i.e., degree $1$ vertices $x$ having $m_{\phi}(x) = 0$).
On the other hand, suppose $x$ is a vertical leaf, and let $e = xy$ be the unique edge with $x \in e$.
Then if $\overline{G}$ is the graph obtained by contracting $e$ to $y$, the induced map
$\overline{G} \to G'$ is still harmonic and non-constant, and $\deg(R_{\overline{G}}) = \deg(R_G)$.
Continuing in this way, we can reduce our claim to the 
already established case where $G$ has no vertical leaves.  
\end{proof}

\begin{remark}
In the classical Riemann-Hurwitz formula from algebraic geometry,
for a non-constant holomorphic map $\phi : X \to X'$ between Riemann surfaces of genus $g$ and $g'$, respectively, one has
\begin{equation*}
2g-2=\deg(\phi)(2g'-2)+\sum_{x \in X}\left( e_{\phi}(x)-1 \right),
\end{equation*}
where $e_{\phi}(x)$ denotes the ramification index of $\phi$ at $x$.
Note that there is no analogue in algebraic geometry of the ``vertical multiplicities'' $v_{\phi}(x)$, and there is an extra factor of $2$ in the 
right-hand side of (\ref{e:RH2}).
Also, note that for Riemann surfaces one has a {\em linear equivalence} $K_{X} \sim \phi^*K_{X'} + R_{X}$
(which is all that can be expected, since there are just canonical {\em divisor classes} on $X$ and $X'$, not canonical
divisors), but in (\ref{e:RiemannHurwitz}) we have an actual equality of divisors.
\end{remark}

\section{Examples}
\label{ExampleSection}

In this section, we give some examples of harmonic and non-harmonic morphisms.

\begin{example}[A harmonic morphism]
The morphism shown in Figure~\ref{Figure1} is harmonic,
with horizontal and vertical multiplicities $m_{\phi}(x)$ and $v_{\phi}(x)$, respectively, labeled next to the
corresponding vertices.
\begin{center}
\begin{figure}[!ht]
\scalebox{0.9}{\includegraphics{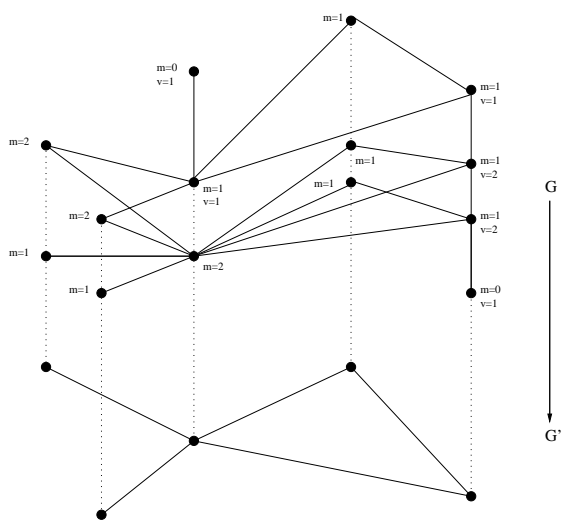}} \caption{A harmonic morphism $\phi : G \to G'$ of degree $3$.}
\label{Figure1}
\end{figure}
\end{center}
\end{example}

\begin{example}[Harmonic morphisms to trees]
\label{HarmonicToTreeExample}
Every graph $G$ admits a non-constant harmonic morphism to a tree.  More precisely, suppose $|V(G)| \geq 2$ and
let $x \in V(G)$ be a vertex of degree $k \geq 1$.  Let $T$ be the graph consisting of two vertices $a,b$ connected
by a single edge $e'$, and let $\phi$ be the morphism sending $x$ to $a$ and every $y \in V(G) \backslash \{ x \}$
to $b$, and sending an edge $e \in E(G)$ to $e'$ if $x \in e$, and to $b$ otherwise.  Then $\phi$ is a harmonic
morphism of degree $k$.
\end{example}

\begin{example}[Automorphisms]
A graph automorphism $\alpha : G \to G$ is a non-degenerate harmonic morphism of degree $1$.
\end{example}

\begin{example}[Coverings]
A morphism $\phi : G \to G'$ is a {\em covering} 
of degree $d\geq 1$ if $\deg(x) = \deg(\phi(x))$ for every $x \in V(G)$ and
$\phi^{-1}(e')$ consists of $d$ disjoint edges for every edge $e' \in E(G')$.
A covering is a harmonic morphism; more precisely, a covering morphism is the same
thing as
a harmonic morphism for which $m_{\phi}(x) = 1$ and
$v_{\phi}(x) = 0$ for all $x \in V(G)$.
\end{example}

\begin{example}[Collapsing]
Let $p \in V(G)$ be a cut vertex, so that $G$ can be partitioned into two subsets $G_1$ and $G_2$
which intersect only at $p$.  We define the {\em collapsing} of $G$ relative to $G_1$
to be the graph $G'$ obtained by contracting all vertices and edges in $G_1$ to $ \{ p \}$.  Let $\phi : G \to G'$
be the morphism which sends $G_1$ to $p$ and is the identity on $G_2$.
Then if $|V(G_2)| > 1$, it is easy to see that $\phi$ is a harmonic morphism of degree 1
(c.f.~\cite[Proposition 4.2]{Urakawa1}).
\end{example}

\begin{example}[Contracting bridges is {\em not} harmonic]
\label{ContractionExample}
Let $e \in E(G)$ be a bridge, 
and let $\overline{G}$ be the graph obtained by contracting $e$.
Then there is an evident {\em contraction morphism} $\rho : G \to \overline{G}$
which is surjective on both vertices and edges.
However, $\rho$ is not in general a harmonic morphism, 
as in Figure~\ref{Figure2}.

\begin{center}
\begin{figure}[!ht]
\scalebox{0.9}{\includegraphics{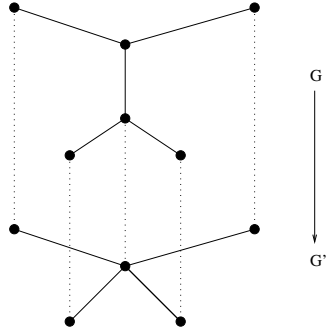}} \caption{A non-harmonic morphism $\rho : G \to G'$.}
\label{Figure2}
\end{figure}
\end{center}

\end{example}

\section{Functorial properties}
\label{FunctorialSection}

In this section, we discuss how harmonic morphisms between graphs induce different kinds of functorial maps between divisor groups, 
Jacobians, and harmonic $1$-forms.

\subsection{Induced maps on Jacobians}

In \S\ref{RiemannHurwitzSection}, we introduced homomorphisms
$\phi_*: \Div(G) \rightarrow \Div(G')$ and
$\phi^*: \Div(G') \rightarrow \Div(G)$ associated to a harmonic morphism $\phi : G \to G'$.
These homomorphisms are related by the following simple formula:

\begin{lem}\label{l:PullPushComposition} Let $\phi: G \rightarrow G'$ be a harmonic morphism, and let $D' \in \Div(G')$.
Then $\phi_*(\phi^*(D'))=\deg(\phi)D'$.
\end{lem}

\begin{proof}
This follows from Lemma~\ref{lem:vertexindependence} and the definitions of $\phi_*$ and $\phi^*$.
\end{proof}

Suppose $\phi: G \rightarrow G'$ is a harmonic morphism and that $f: V(G) \rightarrow A$ and
$f': V(G') \rightarrow A$ are functions, where $A$ is an abelian group.  We define
$\phi_* f : V(G') \to A$ by
\[
\phi_* f(y) := \sum_{\substack{x \in V(G) \\  \phi(x)=y}} m_{\phi}(x) f(x)
\]
and $\phi^* f' : V(G) \to A$ by
\[
\phi^* f' := f' \circ \phi.
\]

\begin{prop}\label{p:PullPushFunct}
Let $\phi: G \rightarrow G'$ be a harmonic morphism, let $f: V(G)
\rightarrow \ZZ$ and $f' \in V(G') \rightarrow \ZZ$.
Then
\begin{equation}\label{e:PullFunctDivisors}
\phi_*(\ddiv(f)) = \ddiv(\phi_*f)
\end{equation}
and
\begin{equation}\label{e:PushFunctDivisors}
\phi^*(\ddiv(f')) = \ddiv(\phi^*f').
\end{equation}
\end{prop}

\begin{proof}
We start by proving (\ref{e:PullFunctDivisors}). We have $$\ddiv(f)=
\sum_{e=xy \in E(G)}(f(x)-f(y))((x)-(y)).$$ By the linearity of $\phi_*$,
we have
\begin{equation}\label{e:PullFunctLeft}\phi_*(\ddiv(f))=
\sum_{e=xy \in E(G)}(f(x)-f(y))(\phi(x)-\phi(y)).
\end{equation}
By the definition of $\phi_*f$, we have
\begin{equation}\label{e:PullFunctRight}
\ddiv(\phi_*f)= \sum_{e'=x'y' \in E(G')} \left( \sum_{x \in V(G),\;
\phi(x)=x'}m_{\phi}(x)f(x) - \sum_{y \in V(G),\;
\phi(y)=y'}m_{\phi}(y)f(y)\right) \left((x')-(y')\right).
\end{equation}
Note that terms in (\ref{e:PullFunctLeft}) corresponding to edges in
$\phi^{-1}(V(G'))$ are zero. Therefore, to derive
(\ref{e:PullFunctDivisors}) from (\ref{e:PullFunctLeft})  and
(\ref{e:PullFunctRight}), it suffices to prove that
$$\sum_{e=xy \in \phi^{-1}(e')}(f(x)-f(y)) = \sum_{x \in V(G),\;
\phi(x)=x'}m_{\phi}(x)f(x) - \sum_{y \in V(G),\;
\phi(y)=y'}m_{\phi}(y)f(y)$$ for every edge $e'=x'y' \in E(G')$. This
last identity holds by the definition of $m_{\phi}$.

We now prove (\ref{e:PushFunctDivisors}). Let $D':=\ddiv(f')$. We have
$D'(y)=\deg(y)f'(y) - \sum_{e=zy \in E(G')}f'(z)$ for every $y \in
V(G')$, so by the definition of $\phi^*$, we have
\begin{equation}\label{e:PushFunctLeft}
(\phi^*D')(x)=m_{\phi}(x)D'(\phi(x))=m_{\phi}(x)\deg(\phi(x))f'(\phi(x))
- m_{\phi}(x)\sum_{e=z\phi(x) \in E(G')}f'(z)
\end{equation}
for every $x \in V(G)$. We now consider $\ddiv(\phi^*f')(x)$. We have
$$\ddiv(\phi^*f')(x)=\ddiv(f' \circ \phi)(x)=\deg(x)f'(\phi(x)) - \sum_{e=xy \in E(G)} f'(\phi(y)).$$
By (\ref{e:DegreeIdentity}), we have $$\deg(x)f'(\phi(x))=
m_{\phi}(x)\deg(\phi(x))f'(\phi(x)) + \sum_{e=xy \in E(G),\;
\phi(y)=\phi(x)}f'(\phi(y)).$$ Therefore
\begin{equation}\label{e:PushFunctRight}
\ddiv(\phi^*f')(x) = m_{\phi}(x)\deg(\phi(x))f'(\phi(x)) - \sum_{e=xy \in E(G),\; \phi(y) \ne \phi(x)}f'(\phi(y)).
\end{equation}
Moreover, for every edge $e'=z\phi(x) \in E(G')$ we have
$$\sum_{e=xy,\; \phi(e)=e'}f'(\phi(y)) = m_{\phi}(x)f'(z),$$
and therefore
$$\sum_{e=xy \in E(G),\; \phi(y) \ne \phi(x)}f'(\phi(y)) = m_{\phi}(x)\sum_{e'=z\phi(x) \in E(G')}f'(z).$$
Thus (\ref{e:PushFunctDivisors}) follows from (\ref{e:PushFunctLeft}) and (\ref{e:PushFunctRight}).
\end{proof}

In particular:

\begin{cor}
\label{cor:principal}
If $\phi : G \to G'$ is a harmonic morphism, then $\phi_*(\Prin(G)) \subseteq \Prin(G')$ and $\phi^*(\Prin(G')) \subseteq \Prin(G)$.
\end{cor}

As a consequence of Corollary~\ref{cor:principal},
we see that $\phi$ induces group homomorphisms (which we continue to denote by $\phi_*,\phi^*$)
\[
\phi_* : \Jac(G) \to \Jac(G'), \,
\phi^* : \Jac(G') \to \Jac(G).
\]
It is straightforward to check that if $\psi : G \to G'$ and $\phi : G' \to G''$ are harmonic morphisms and
$D \in \Div(G), D'' \in \Div(G'')$, then $\phi \circ \psi : G \to G''$ is harmonic, and we have
$(\phi \circ \psi)_*(D) = \phi_*(\psi_*(D))$ and
$(\phi \circ \psi)^*(D'') = \psi^*(\phi^*(D''))$.
Therefore we obtain two different functors from the category of graphs (together with harmonic morphisms between them)
to the category of abelian groups: a covariant ``Albanese'' functor $(G \mapsto \Jac(G), \phi \mapsto \phi_*)$ and a
contravariant ``Picard'' functor $(G \mapsto \Jac(G), \phi \mapsto \phi^*)$.
(This terminology comes from the corresponding notions in algebraic geometry.)

\medskip

\begin{remark}
\label{AlbaneseRemark}
A more conceptual definition of the Albanese functor $\phi_*$ is as follows.  Choose a base vertex $x_0 \in G$, and let $S = S_{x_0} :
G \to \Jac(G)$ and $S' = S_{\phi(x_0)} : G' \to \Jac(G')$
denote the corresponding Abel-Jacobi maps.  Since $S' : G' \to \Jac(G')$ is a harmonic function, it follows from 
Proposition~\ref{p:conformity} that the pullback $S' \circ \phi$ is a harmonic map from $G$ to $\Jac(G')$. 
As $S' \circ \phi$ sends $x_0$ to $0$, it follows from Lemma~\ref{l:AbelJacobiUnivProperty}
that there is a unique homomorphism $\psi : \Jac(G) \to \Jac(G')$
such that $S' \circ \phi = \psi \circ S$.  
From the uniqueness of $\psi$, it follows easily that $\psi = \phi_*$.

In particular, a harmonic morphism $\phi : G \to G'$ gives rise to a commutative diagram
\[
\begin{CD}
G     @>{\phi}>>     G'     \\
@V{S}VV    @VV{S'}V     \\
\Jac(G)    @>{\phi_*}>>   \Jac(G')   \\
\end{CD}
\]
\end{remark}

\medskip

As an application of Corollary~\ref{cor:principal}, we have the following result:

\begin{cor}
\label{cor:grdtransfer}
If $\phi : G \rightarrow G'$ is a non-constant harmonic morphism,
then for every $D \in \Div(G)$ we have $r_{G'}(\phi_*(D)) \geq
r_G(D)$.
\end{cor}

\begin{proof}
By Lemma~\ref{lem:constantorsurj}, $\phi$ is surjective on vertices.
Let $D' := \phi_*(D)$, and let $k$ be a nonnegative integer.  For
every effective divisor $E' \in \Div(G')$ of degree $k$, we can
choose $E \in \Div(G)$ such that $\phi_*(E) = E'$. If $r_G(D) \geq
k$, then by definition $D-E = F + P$ with $F$ effective and $P$
principal.  Applying $\phi_*$ and using the fact that $\phi_*(P) \in
\Prin(G')$, we see that $D' - E'$ is equivalent to the effective
divisor $\phi_*(F)$, and therefore $r_{G'}(D') \geq k$ as well.
\end{proof}

\medskip

We now investigate some useful general properties of the induced maps
$\phi_*$ and $\phi^*$ on Jacobians.  In the classical algebraic geometry setting,
$\phi_*$ is always surjective but $\phi^*$ is sometimes injective and sometimes not.
More precisely, $\phi^* : \Jac(X') \to \Jac(X)$ is injective if and only if
$\phi : X \to X'$ has a nontrivial unramified abelian subcover.  The situation for graphs
is simpler, since as we will now show, $\phi_*$ is always surjective and $\phi^*$ is always injective.
The surjectivity of $\phi_*$ is easy:

\begin{lem}
\label{lem:surjlem}
Let $\phi: G \rightarrow G'$ be a non-constant harmonic morphism.  Then
$\phi_* : \Jac(G) \to \Jac(G')$ is surjective.
\end{lem}

\begin{proof}
It follows from Lemma~\ref{lem:constantorsurj} and the linearity
of $\phi_*$ that  $\phi_*$ is a surjective map from $\Div(G)$ to
$\Div(G')$, which implies surjectivity on the level of Jacobians.
\end{proof}

The injectivity of $\phi^*$ is much more subtle (as one would expect, since the
analogous statement for Riemann surfaces is false):

\begin{thm}
\label{thm:injthm}
Let $\phi: G \rightarrow G'$ be a non-constant harmonic morphism.  Then
$\phi^* : \Jac(G') \to \Jac(G)$ is injective.
\end{thm}

\begin{proof}
We first set the following notation.  For a function $f: V(G) \to \ZZ$, let $\max(f) =
\max_{x \in V(G)}f(x)$, let $\min(f) = \min_{x \in V(G)}f(x)$, and
let $s(f) = \max(f) - \min(f)$. Let $M(f)=\{x \in V(G) \; | f(x)=
\max(f)\}$, and let $m(f)=\{x \in V(G) \; | f(x)= \min(f)\}$.

It suffices to show that $D' \in \Prin(G')$ for every $D'
\in \Div(G')$ such that $\phi^*(D') \in \Prin(G)$. Suppose for the sake of
contradiction that there exists a divisor $D' \in \Div(G') \setminus
\Prin(G')$ such that $\phi^*(D') = \ddiv(f)$ for some $f: V(G) \to
\ZZ$. Choose such a $D'$ for which $s(f)$ is minimized, and
subject to this condition such that $|M(f)|$ is minimized. Let $D :=
\phi^*(D') = \ddiv(f)$.  Clearly $s(f) \neq 0$, as otherwise $D = 0$, and
therefore $D'=0 \in \Prin(G')$, a contradiction. Therefore there exists 
a vertex $x_0 \in M(f)$ with a neighbor in $V(G) \setminus M(f)$. 
For every $x \in M(f)$, one has
$$D(x) = \ddiv(f)(x) \geq |\{e \in E(G) \; | \; e=xy, y \in V(G) \setminus
M(f)\}|.$$ It follows that $D(x) \geq 0$ for every $x \in M(f)$, and
that $D(x_0)>0$. Similarly,  for every $x \in m(f)$ one has either
$D(x) < 0$, or else $D(x)=0$ and all the neighbors of $x$ belong to $m(f)$. Let
$X=\phi^{-1}(\phi(x_0)) \cap m(f)$.
Since $D(x_0) > 0$, we have $D'(\phi(x_0))>0$ as well, so by the definition of $\phi^*$
it follows that $D(x)>0$ for every $x \in \phi^{-1}(\phi(x_0))$ with $m_\phi(x) > 0$.
Therefore $X$ consists entirely of vertices
$x \in \phi^{-1}(\phi(x_0))$ with $m_{\phi}(x) = 0$ and $D(x) = 0$.
But then all the neighbors of vertices in $X$ belong to
$\phi^{-1}(\phi(x_0))$, and thus by the above must belong to $X$. 
Since $G$ is connected, it follows that $X$ is empty, i.e.,
\begin{equation}\label{e:nonMinimum}
\phi^{-1}(\phi(x_0)) \cap m(f) = \emptyset.
\end{equation}

Let $\chi: V(G') \to \ZZ$ be the characteristic function of
$\{\phi(x_0)\}$,
and let $D'' = D' - \ddiv(\chi)$. We claim that $D''$ contradicts the choice of $D'$.
Clearly, $D'' \in \Div(G') \setminus \Prin(G')$. By
Proposition~\ref{p:PullPushFunct}, we have $$\phi^*(D'') = \phi^*(D') -
\phi^*(\ddiv(\chi))= \ddiv(f)-\ddiv(\phi^*\chi)=\ddiv(f -
\chi \circ \phi).$$ Let $D^\star = \phi^*(D'')$ and let $f^\star = f -
\chi \circ \phi$. We have $f^\star(x) =  f(x)-1$ for every $x \in
\phi^{-1}(\phi(x_0))$, and $f^\star(x) = f(x)$ otherwise.
By (\ref{e:nonMinimum}), we have $\min(f) = \min(f^\star)$, and
clearly $\max(f) \geq \max(f^\star)$. Therefore $s(f) \geq
s(f^\star)$. Moreover, either $s(f) > s(f^\star)$ or $\max(f) =
\max(f^\star)$. In the second case, we have $M(f^\star) \subseteq
M(f) \setminus \{x_0\}$, and thus $|M(f)|
> |M(f^\star)|$. It follows that $D''$ contradicts the choice of $D'$, as claimed.
\end{proof}

\subsection{Eulerian cuts and harmonic morphisms}

Let $\kappa_G = |\Jac(G)|$ denote the number of spanning trees in a
graph $G$. From either Lemma~\ref{lem:surjlem} or Theorem~\ref{thm:injthm}, we
immediately deduce the following corollary:

\begin{cor}
\label{cor:spanningtrees}
If there exists a non-constant harmonic
morphism from $G$ to $G'$, then $\kappa_{G'}$ divides $\kappa_G$.
\end{cor}

Define an {\em Eulerian cut} in a graph $G$ to be a non-empty cut 
which is also an even subgraph of $G$; equivalently,
an Eulerian cut is a partition of $V(G)$ into non-empty disjoint subsets
$X$ and $X'$ in such a way that there are an even number of edges connecting each vertex in $X$ (resp. $X'$) to
vertices in $X'$ (resp. $X$).  According to a theorem of Chen 
\cite{Chen} (see also \cite[Proposition~35.2]{BiggsAPTG}),
$G$ has an Eulerian cut if and only if $\kappa_G$ is even.
From Corollary~\ref{cor:spanningtrees}, it therefore follows that
if $G'$ has an Eulerian cut and there exists a non-constant harmonic
morphism from $G$ to $G'$, then $G$ has an Eulerian cut as well.
We can strengthen this observation using the following result, which
characterizes the existence of Eulerian cuts in terms of 
non-constant harmonic maps from $G$ to the graph $B_2$ consisting of $2$ vertices connected by $2$ edges:

\begin{thm}
\label{thm:paritythm}
Let $G$ be a graph.  Then the following are equivalent:
\begin{enumerate}
\item $G$ has an Eulerian cut.
\item There is a non-constant harmonic morphism from $G$ to $B_2$.
\item $\kappa_G$ is even.
\end{enumerate}
\end{thm}

\begin{proof}
Although the equivalence $(1) \Leftrightarrow (3)$ is just Chen's theorem, for the reader's convenience we will 
provide a self-contained proof of this result.  Our proof of $(3) \Rightarrow (1)$ is borrowed from
the unpublished manuscript \cite{Eppstein}.

${\bf (1) \Rightarrow (2):}$ 
Suppose that $G$ admits an Eulerian cut $S$. 
We claim that there exists a partition of $S$ into non-empty disjoint subsets 
$S_1, S_2 \subseteq E(G)$ such that
$$|\{e \in S_1 \; | \; x \in e \}| = |\{e \in S_2 \; | \; x \in e \}|$$
for every $x \in V(G)$. 
Indeed, it is well-known (see \cite[\S{I.1},Theorem~1]{Bollobas}) 
that the edge set of any Eulerian graph can be decomposed into edge-disjoint cycles. 
Since the graph $G[S]$ with vertex set $V(G)$ and edge set $S$ is
Eulerian and bipartite, it follows that $G[S]$ decomposes into
edge-disjoint {\em even} cycles. It is trivial to obtain the required
partition for an even cycle. By composing the resulting partitions
of even cycles, one then obtains the required partition $(S_1,S_2)$ of
$S$.

We now construct a non-constant harmonic morphism $\phi: G \to B_2$ as
follows. Let the vertices of $B_2$ be labeled $x$ and $y$, and let
the edges of $B_2$ be labeled $e_1$ and $e_2$. Let $X \subseteq
V(G)$ be one of the sides of $S$. For $z \in V(G)$, let $\phi(z)=x$
if $x \in X$, and let $\phi(z)=y$ otherwise. For $i \in \{1,2\}$
and $e \in S_i$, let $\phi(e) = e_i$. Finally, if $e = z_1z_2 \in
E(G) \setminus S$, we set $\phi(e)=\phi(z_1)=\phi(z_2)$. It
follows from the definition of $S_1$ and $S_2$ that $\phi$ is a
non-constant harmonic morphism.

${\bf (2) \Rightarrow (3):}$ 
We have $\kappa_{B_2}=2$. Therefore, if $G$ admits a non-constant
harmonic morphism to $B_2$, then $\kappa_G$ is even by
Corollary~\ref{cor:spanningtrees}.

${\bf (3) \Rightarrow (1):}$ 
(c.f.~\cite[Proof of Theorem~6]{Eppstein})
Let $\Lambda(G) = H^1(G,\ZZ) \subset H^1(G,\RR)$ denote 
the {\em lattice of integral flows} on $G$, and let $\Lambda^{\#}(G)$ be the 
lattice dual to $\Lambda(G)$ under the standard Euclidean inner product 
$\langle \, , \, \rangle$ on $C^1(G,\RR) \supseteq H^1(G,\RR)$.
Explicitly, we have 
\[
\Lambda^{\#} = \{ \omega \in H^1(G,\RR) \; | \; \langle \omega, \omega' \rangle \in \ZZ
\textrm{ for all } \omega' \in H^1(G,\ZZ) \}.
\]
By \cite{BDN} (see also \cite[\S{29}]{BiggsAPTG}), there is 
a canonical isomorphism $\Jac(G) \cong \Lambda^{\#}(G) / \Lambda(G)$.

Suppose that $\kappa_G$ is even.  Then $\Jac(G)$ has an element of order 2, so there is a flow
$\omega \in \Lambda^{\#}(G)$ such that $\omega \not\in \Lambda(G)$ but $2\omega \in \Lambda(G)$.
Thus the value of $\omega$ along each edge of $G$ is a half-integer, 
and the set $S$ of edges along which $\omega$ is non-integral is non-empty.
Since $\delta(\omega) = 0$, it follows that every vertex in $S$ has even degree.  So it
suffices to prove that $S$ is a cut.
To see this, choose a vertex $x \in V(G)$, and partition $V(G)$ into disjoint subsets $A$ and $B$
as follows: a vertex $y \in V(G)$ belongs to $A$ (resp. $B$) iff it can be connected to $x$ by a path containing
an {\em odd} (resp. {\em even}) number of edges in $S$.
Since $G$ is connected, $A \cup B = V(G)$.
Furthermore, we have $A \cap B = \emptyset$, because otherwise there would be a cycle $C$ in
$G$ containing an odd number of edges of $S$, and therefore $\langle \omega,\chi_C \rangle \not\in \ZZ$,
contradicting the fact that $\omega \in \Lambda^{\#}(G)$.  
Finally, to see that $S$ is indeed a cut, note that each edge $e \in S$ connects a vertex in $A$ to a vertex in $B$
(since $e$ itself is a path with one edge in $S$), and an edge $e' \not\in S$ cannot connect a vertex in $A$
to a vertex in $B$ (since $e'$ is a path with no edges in $S$).  Thus $S$ is precisely the cut consisting of
all edges connecting $A$ to $B$.
\end{proof}

\begin{remark}
Here is a more direct argument for proving $(1) \Rightarrow (3)$ which makes use of Theorem~\ref{t:OrderDivisor}.
Let $S$ be an Eulerian cut in $G$ separating the subsets $X,Y \subset V(G)$. It is easy to see that 
there exists an ordering $x_1,\ldots,x_k$ of $X$ such that for every $i \in \{1, \ldots, k\}$ either 
$\outdeg_X(x_i)>0$ or $x_jx_i \in E(G)$ for some $j<i$. Similarly, there exists an ordering $y_1,\ldots,y_\ell$ of $Y$ such that 
for every $i\in \{1, \ldots, l\}$ either $\outdeg_Y(y_i)>0$ or $y_jy_i \in E(G)$ for some $j<i$. 
Define a divisor $D \in \Div^0(G)$ by setting $D(x) := \frac{1}{2} \outdeg_X(x)$ for $x \in X$, and $D(y) := -\frac{1}{2} \outdeg_Y(y)$ for $y \in Y$.  
Then $2D = \ddiv(\chi_X) \sim 0$. However, using Theorem~\ref{t:OrderDivisor}, we see that $D$ itself is not equivalent to $0$, since $D \leq \nu$, 
where $\nu$ is the divisor associated to the linear order $y_1 < \cdots < y_\ell < x_1 < \cdots < x_k$ on $V(G)$.
Thus $D$ corresponds to an element of order $2$ in $\Jac(G)$, and in particular $\kappa_G = |\Jac(G)|$ is even.
\end{remark}

\medskip

\subsection{Induced maps on harmonic $1$-forms}

We now turn to a discussion of harmonic $1$-forms and the maps induced on them
by a harmonic morphism.

\medskip

We begin with some notation and terminology.
Let $\vec{E}(G)$ denote the set of directed edges of $G$.  
For $e \in \vec{E}(G)$, we
let $o(e),t(e)$ denote the {\em origin} and {\em terminus} of $e$, respectively.
We denote by $\overline{e}$ the directed edge representing the same undirected edge as $e$, but
with the opposite orientation.
From the definition of a morphism of graphs, it follows easily that 
a morphism $\phi: G \to G'$ induces a natural map from $\vec{E}(G)$ to
$\vec{E}(G') \cup V(G')$.

\medskip

Let $A$ be an abelian group, and let $C^1(G,A)$ denote the space of $1$-cochains on $G$ with
values in $A$, i.e., functions $\omega : \vec{E}(G) \to A$ with the property that 
$\omega(e) = -\omega(\overline{e})$ for all $e \in \vec{E}(G)$.
As usual, we also let $C^0(G,A)$ denote the space of all functions $f : V(G) \to A$.
We define the {\em coboundary operator} $\delta : C^1(G,A) \to C^0(G,A)$ by the formula
\begin{equation}
\label{e:coboundary}
\delta(\omega)(x) := \sum_{\substack{e \in \vec{E}(G) \\ t(e) = x}} \omega(e).
\end{equation}

An {\em $A$-flow} (or simply a {\em flow} if $A = \RR$) on $G$ is a $1$-cochain $\omega \in C^1(G,A)$
such that $\delta(\omega) = 0$.  We denote by $H^1(G,A)$ the space of $A$-flows on $G$.  When $A=\RR$,
we will also refer to $\H^1(G) := H^1(G,\RR)$ as the space of {\em harmonic $1$-forms} on $G$;
it is analogous to the space $\Omega^1(X)$ of holomorphic $1$-forms on a Riemann surface $X$.  For example, it is
well-known that $\dim_{\RR} \H^1(G) = g$ (just as $\dim_{\CC} \Omega^1(X) = g$ in the Riemann surface case).

\medskip

Suppose $\phi: G \rightarrow G'$ is a harmonic morphism and that 
$\omega \in C^1(G,A),\omega' \in C^1(G',A)$ are $1$-cochains.
We define the {\em pullback}
$\phi^* \omega' \in C^1(G,A)$ by
\[
(\phi^* \omega')(e) := \left\{
\begin{array}{ll}
\omega'(\phi(e)) & \textrm{ if } \phi(e) \in \vec{E}(G') \\
0 & \textrm{ otherwise}
\end{array} \right.
\]
and the {\em push-forward} (or {\em trace}) 
$\phi_* \omega \in C^1(G',A)$ by
\[
\phi_* \omega(e') := \sum_{\substack{e \in \vec{E}(G) \\  \phi(e)=e'}} \omega(e).
\]

\begin{prop}\label{p:PullPushDiff}
Let $\phi: G \rightarrow G'$ be a harmonic morphism
and let $\omega \in \H^1(G), \omega' \in \H^1(G')$ be harmonic $1$-forms.
Then: 
\begin{enumerate}
\item $\phi^*\omega' \in \H^1(G)$.
\item $\phi_*\omega \in \H^1(G')$.
\end{enumerate}
\end{prop}

\begin{proof}
To establish (1), we follow \cite[Proof of Theorem~2.13]{Urakawa1}.
For every $x \in V(G)$, we have
\begin{align*}
\sum_{\substack{e \in \vec{E}(G), t(e) = x \\ \phi(e) \in \vec{E}(G')}} \omega'(\phi(e))
= \sum_{\substack{e' \in \vec{E}(G') \\ t(e') = \phi(x)}} \sum_{\substack{e \in \vec{E}(G), x \in e \\ \phi(e) = e'}} \omega'(\phi(e)) 
= m_{\phi}(x) \sum_{\substack{e' \in \vec{E}(G') \\ t(e') = \phi(x)}} \omega'(e'). \\
\end{align*}
Since $\delta(\omega') = 0$ and $(\phi^* \omega')(e) = 0$ for all vertical edges $e \in \vec{E}(G)$, for all $x \in V(G)$ we have
\begin{align*}
\delta(\phi^* \omega')(x) &= 
\sum_{\substack{e \in \vec{E}(G) \\ t(e) = x}} (\phi^*\omega')(e) 
= \sum_{\substack{e \in \vec{E}(G), t(e) = x \\ \phi(e) \in \vec{E}(G')}} \omega'(\phi(e)) \\
&= m_{\phi}(x) \sum_{\substack{e' \in \vec{E}(G') \\ t(e') = \phi(x)}} \omega'(e') 
= m_{\phi}(x) \delta(\omega')(\phi(x)) \\
&= 0, 
\end{align*}
which proves (1).

For (2), note that for every $y \in V(G')$, we have
\begin{equation}
\label{e:deltacalc}
\sum_{\substack{e' \in \vec{E}(G') \\ t(e') = y}} 
\sum_{\substack{e \in \vec{E}(G) \\ \phi(e) = e'}} \omega(e) 
= \sum_{\substack{x \in V(G) \\ \phi(x) = y}} 
\sum_{\substack{e \in \vec{E}(G) \\ t(e) = x}} \omega(e),
\end{equation}
since each vertical edge in $E(G)$ gets counted twice in the sum on the right-hand side of (\ref{e:deltacalc}), once with each orientation,
and therefore the net contribution to the sum from such an edge is zero.
Therefore
\begin{align*}
\delta(\phi_* \omega)(y) &= 
\sum_{\substack{e' \in \vec{E}(G') \\ t(e) = y}} (\phi_*\omega)(e') 
= \sum_{\substack{e' \in \vec{E}(G') \\ t(e') = y}} 
\sum_{\substack{e \in \vec{E}(G) \\ \phi(e) = e'}} \omega(e) \\
&= \sum_{\substack{x \in V(G) \\ \phi(x) = y}} 
\sum_{\substack{e \in \vec{E}(G) \\ t(e) = x}} \omega(e) \qquad\qquad\text{by (\ref{e:deltacalc})} \\
&= \sum_{\substack{ x \in V(G) \\ \phi(x) = y}} \delta(\omega)(x) = 0,
\end{align*}
proving (2).
\end{proof}

As a consequence of Proposition~\ref{p:PullPushDiff},
we see that $\phi$ induces linear transformations (which we continue to denote by $\phi^*,\phi_*$)
\[
\phi^* : \H^1(G') \to \H^1(G), \,
\phi_* : \H^1(G) \to \H^1(G').
\]
It is straightforward to check that the association $(G',\phi) \mapsto (\H^1(G'),\phi^*)$ 
(resp. $(G,\phi) \mapsto (\H^1(G'),\phi_*)$) is a contravariant (resp. covariant) functor from 
the category of graphs (together with harmonic morphisms between them)
to the category of vector spaces.

It follows easily from the definitions that 
\begin{equation}
\label{e:PullPushComposition2}
\phi_*\phi^*(\omega') = \deg(\phi) \omega'
\end{equation}
for all $\omega' \in \H^1(G')$ (compare with Lemma~\ref{l:PullPushComposition}).
As a consequence, we obtain the following result,
which provides another way to see that if 
$\phi$ is a non-constant harmonic morphism from a graph of genus $g$ to a graph of genus $g'$, then $g' \leq g$
(c.f.~Theorem~\ref{t:RiemannHurwitz}):

\begin{cor}
\label{cor:injsurj}
If $\phi : G \to G'$ is a non-constant harmonic morphism, then $\phi^* : \H^1(G') \to \H^1(G)$ is injective
and $\phi_* : \H^1(G) \to \H^1(G')$ is surjective. 
\end{cor}

\begin{proof}
Both the injectivity of $\phi^*$ and the surjectivity of $\phi_*$ follow easily from (\ref{e:PullPushComposition2}).
However, one can also prove the injectivity of $\phi^*$ directly (c.f.~\cite[Proof of Theorem 2.13]{Urakawa1}): if 
$\phi^*(\omega') = 0$, then $\omega'(\phi(e)) = 0$ for all horizontal edges $e \in E(G)$, and
since $\phi$ maps the set of horizontal edges of $G$ surjectively onto $E(G')$, it follows that 
$\omega' = 0$.  
\end{proof}

By functoriality, an automorphism $\alpha$ of a graph $G$ induces an automorphism 
$\alpha^*$ of the vector space $\H^1(G)$.
For later use, we note the following property of the corresponding map $\Aut(G) \to \Aut(\H^1(G))$:

\begin{prop}
\label{prop:autinj}
If $G$ is a $2$-edge-connected graph of genus at least $2$, 
then the natural map from $\Aut(G)$ to $\Aut(\H^1(G))$ is injective.
\end{prop}

\begin{proof}
Let $\beta,\beta' \in \Aut(G)$.
By considering the automorphism $\alpha := \beta' \beta^{-1}$, it suffices to prove that 
if $\alpha^*$ is the identity map on $\H^1(G)$, then $\alpha$ is the identity map on $G$.
So suppose $\alpha^* = {\rm Id}$.  Then every directed cycle in $G$ is mapped onto itself.
Let $C$ be an (undirected) simple cycle in $G$ (i.e., a cycle with no repeated vertices), 
let $x \in C$ be a vertex of degree at least 3, and let $x' = \alpha(x)$.
Let $e \in C$ be the directed edge with $o(e) = x$, let $e' = \alpha(e) \in C$, 
and let $e'' \in \vec{E}(G)$ be a directed edge with $o(e'') = x$ and $e'' \not\in C$.  
Since $G$ is $2$-edge-connected, $e''$ belongs to a simple cycle $C''$, and
we can choose $C''$ so that either $V(C) \cap V(C'') = \{ x \}$,
or else so that $E(C) \cap E(C'')$ is a path in $C$ containing $e$.

{\bf Case I:} $V(C) \cap V(C'') = \{ x \}$.

In this case, $x' \in V(C) \cap V(C'') = \{ x \}$ so $\alpha(x) = x$.  
But then $\alpha(e) = e$, since $\alpha^*$ preserves {\em directed} cycles of $G$.
From this it follows easily that $\alpha$ is the identity map on $C$.

{\bf Case II:} $E(C) \cap E(C'')$ is a path in $C$ containing $e$.

In this case, we must also have $e' \in C''$.  Suppose $e' \neq e$.  
Then as $\alpha(e'') \not\in C$, the cycle $C''$ can be directed so that it consists 
of the unique path in $C$ from $x$ to $x'$ followed by the unique path in $C'' \backslash C$ from $x'$ to $x$.
But then $\alpha$ restricted to $C''$ is orientation-reversing, a contradiction.  
We conclude that $\alpha(e) = e$, and hence $\alpha$ is the identity map on $C$ in this case as well.

It follows that the restriction of $\alpha$ to every simple cycle $C$ of $G$ is the identity map.
Since $G$ is $2$-edge-connected, this implies that $\alpha$ is the identity map on all of $G$.
\end{proof}

\begin{remark}
Proposition~\ref{prop:autinj} is the analogue of the fact from algebraic geometry that if $X$ is a Riemann surface of genus at least
$2$, then the natural map from $\Aut(X)$ to $\Aut(\Omega^1(X))$ is injective.
\end{remark}

As a consequence of Proposition~\ref{prop:autinj}, we obtain the following non-trivial restriction
on the automorphism group of a $2$-edge-connected graph of genus at least 2:

\begin{cor}
\label{cor:autbound}
If $G$ is a $2$-edge-connected graph of genus $g \geq 2$, 
then the group $\Aut(G)$ is isomorphic to a subgroup of the group $\GL(g,\ZZ)$ of invertible $g \times g$ matrices 
with coefficients in $\ZZ$.  
\end{cor}

\begin{proof}
Since $\Aut(G)$ acts faithfully on the $g$-dimensional vector space $\H^1(G,\RR)$ and preserves the lattice $\H^1(G,\ZZ)$, the result follows. 
\end{proof}

\begin{remark}
By a theorem of Minkowski, the torsion group $\GL(n,\ZZ)_{\rm tors}$ of $\GL(n,\ZZ)$ is finite for all $n \geq 1$, and every prime divisor $p$ of $|\GL(n,\ZZ)_{\rm tors}|$ satisfies $p \leq n+1$.
In particular, Corollary~\ref{cor:autbound} implies that if a $2$-edge-connected graph $G$ of genus $g \geq 2$ has an automorphism of prime order $p$ then $p \leq g+1$.
This bound is sharp, since the graph $B_{n+1}$ consisting of $2$ vertices joined by $n+1$ edges has genus $n$ and $|\Aut(B_{n+1})| = 2(n+1)!$.
\end{remark}

\section{Hyperelliptic graphs}
\label{HyperellipticSection}

\subsection{Definition and basic properties}

We say that a graph $G$ is \emph{hyperelliptic} if there exists a
divisor $D \in \Div(G)$ such that $\deg(D) = 2$ and $r(D)=1$. By
Riemann-Roch for graphs, if $G$ is hyperelliptic then $g(G) \geq 2$,
and by Clifford's theorem for graphs, if $g(G) \geq 2$ and $\deg(D)
= 2$, then $r(D) = 1$ if and only if $r(D) \geq 1$.

\begin{example}
Every graph of genus $2$ is hyperelliptic.
Indeed, if $g(G) = 2$, then by Riemann-Roch for graphs,
the canonical divisor $K_G$ has $\deg(K_G) = 2$
and $r(K_G) = 1$.
\end{example}

\begin{example}\label{ex:Banana} Let the graph $G=B(l_1,l_2, \ldots, l_n)$ consist
of two vertices $x$ and $y$ and $n \geq 3$ internally disjoint paths
joining $x$ to $y$ with lengths $l_1,l_2, \ldots, l_{n}$.
Then $G$ is hyperelliptic. More specifically, we claim that $r((x)+(y)) = 1$. To prove
this, it suffices to show that $|(x) + (y) -(z)|\neq \emptyset$
for every $z \in V(G)$. Consider one of the paths joining $x$ and
$y$, and let $x, z_1, z_2, \ldots, z_{l-1}, y$ be the vertices of
this path in order. Then $(x)+(y) \sim (z_i) + (z_{l-i})$, and
therefore $|(x) + (y) -(z_i)|\neq \emptyset$ for every $1 \le i \le
l-1$. Thus $r((x)+(y)) = 1$, and our claim follows.
\end{example}

Although the graph $G=B(1,1,\ldots,1)$ has edge connectivity equal to $|E(G)|$, which
can be arbitrarily large, the following result shows that every other hyperelliptic graph has edge connectivity
at most 2:

\begin{lem}
\label{lem:hyperelledgeconn}
If $G$ is a hyperelliptic graph, then either $|V(G)| = 2$ (so that $G$ is isomorphic to a graph of
the form $B(1,1,\ldots,1)$) or $G$ has edge connectivity at most 2.
\end{lem}

\begin{proof}
Let $D = (x) + (x')$ be an effective divisor of degree $2$ on $G$ with $r(D) = 1$.
If $|V(G)| > 2$, choose a vertex $y \in V(G)$ with $y \not\in \{ x,x' \}$.
Since $r(D) = 1$, there exists $y' \in V(G)$ such that $(x) + (x') \sim (y) + (y')$,
and therefore the map $S^{(2)} : \Div_{+}^{2}(G) \to \Jac(G)$ is not injective.
By Theorem~\ref{t:GraphAbelJacobi}, it follows that $G$ is not $3$-edge-connected.
\end{proof}

A classical result from algebraic geometry asserts that if $X$ is a hyperelliptic Riemann surface and $\phi : X \to X'$ is a
non-constant holomorphic map with $g(X') \geq 2$, then $X'$ is also hyperelliptic.
Using Corollary~\ref{cor:grdtransfer}, we obtain the following analogous result for graphs:

\begin{cor}
If $G$ is hyperelliptic and $\phi : G \to G'$ is a non-constant harmonic
morphism onto a graph $G'$ with $g(G') \geq 2$, then $G'$ is hyperelliptic as well.
\end{cor}

\medskip

As in classical algebraic geometry, we can also show in the
graph-theoretic setting that there is at most one complete linear
system $|D|$ of degree 2 on a graph $G$ for which $r(D) = 1$:

\begin{prop}
\label{prop:uniqueg12}
If $D,D'$ are degree 2 divisors on $G$ with $r(D) = r(D') = 1$, then $D \sim D'$.
\end{prop}

\begin{proof}
We may assume that $g := g(G) \geq 2$. Consider the divisor $E := D
+ (g-2)D'$ of degree $2g-2$ on $G$.  By Lemma~\ref{SubAdditivityLemma}, we have
$r(E) \geq g-1$. By Riemann-Roch for graphs, we have $r(K_G-E) \geq
0$; since $\deg(K_G-E) = 0$, it follows that $K_G \sim E$.  Applying
the same reasoning to $E' := (g-1)D'$, we see that $K_G \sim E'$, and
therefore $D \sim D'$ as desired.
\end{proof}

\medskip

\subsection{Hyperelliptic graphs, involutions, and harmonic morphisms}

Our next goal is to obtain a graph-theoretic analogue of the
well-known result from algebraic geometry that the following are
equivalent for a Riemann surface $X$ of genus at least $2$: 
(i) $X$ is  hyperelliptic; 
(ii) $X$ admits a non-constant holomorphic map of degree 2 onto the Riemann sphere; and 
(iii) there is an involution $\iota : X \to X$
whose quotient is isomorphic to the Riemann sphere.
We begin by discussing quotients in the category of graphs
(together with morphisms between them).

\medskip

Let $H$ be a finite group acting on a graph $G$, i.e., suppose we are given a homomorphism $H \to \Aut(G)$.  We write $h \cdot x$ for
the action of an element $h \in H$ on an element $x$ of $V(G) \cup E(G)$.
We define the quotient graph $G/H$, together with a canonical morphism $\pi_H : G \to G/H$, as follows.

\medskip

For $x,y \in V(G) \cup E(G)$, let $x \sim_{H} y$ if there exists an
element $h \in H$ such that $h \cdot x = y$. Then $\sim_{H}$ is an
equivalence relation  on $V(G) \cup E(G)$. The {\em quotient graph}
$G / H$ is constructed as follows. The vertices of $G/H$ are the
equivalence classes of $V(G)$ with respect to $\sim_{H}$. The edges
of $G/H$ correspond to those equivalence classes of $E(G)$ with
respect to $\sim_{H}$ which consist of edges whose ends are
inequivalent.  It is readily verified that $G/H$ is a graph in our
sense of the word (i.e., a connected multigraph with no loop edges).  
The \emph{quotient morphism} $\pi_H : G \rightarrow G/H$
maps every vertex of $G$ to its equivalence class, every edge of $G$
whose ends are inequivalent to the edge of $G/H$ corresponding to
its equivalence class, and every edge of $G$ with equivalent ends to
the equivalence class of its ends.  It is straightforward to check
that $\pi_H$ is a surjective morphism of graphs (though not necessarily
a harmonic morphism), and by construction
we have $\pi_H(h\cdot x) = \pi_H(x)$ for all $h \in H$ and all $x
\in V(G)\cup E(G)$.  In fact, the morphism $\pi_H : G \to G/H$ has the
following universal property: if $\pi' : G \to G/H$ is any morphism
of graphs for which $\pi'(h\cdot x) = \pi'(x)$ for all $h \in H$ and
all $x \in V(G) \cup E(G)$, then there exists a unique morphism $\psi : G/H
\to G'$ such that $\pi' = \psi \circ \pi_H$. This universal property
uniquely characterizes $G/H$ up to isomorphism.

\medskip

If $H = \langle \phi \rangle$ is a cyclic subgroup of $\Aut(G)$, we
will often write $G/{\phi}$ instead of $G/H$ and $\phi^{\sim}$ instead
of $\pi_{\phi}$.

\medskip

An automorphism $\iota$ of a graph $G$ is called an
\emph{involution} if $\iota \circ \iota$ is the identity
automorphism. We say that an involution $\iota$ is \emph{mixing} if
for every edge $e=xy \in E(G)$ such that $\iota(e) = e$ we have
$\iota(x) = y$.  Equivalently, $\iota$ is mixing if and only if it does not fix 
any directed edge of $G$.
The following lemma shows that if $|V(G)|>2$, 
there is a one-to-one correspondence between mixing involutions of $G$ 
and non-degenerate harmonic morphisms of degree two from $G$ to a graph $G'$.


\begin{lem}\label{l:Involution} Let $G,G'$ be graphs, and let $\phi: G \rightarrow G'$ be
a non-degenerate harmonic morphism of degree $2$.  Then there
is a mixing involution $\iota$ of $G$ for which
$\phi=\iota^{\sim}$. Conversely, let $|V(G)|>2$ and let $\iota: G
\rightarrow G$ be a mixing involution.  Then $\iota^{\sim}$ is a
non-degenerate harmonic morphism of degree two.
\end{lem}

\begin{proof} Let $\phi: G \rightarrow G'$ be a non-degenerate
harmonic morphism of degree $2$. For $x \in V(G)$, if there
exists $y \neq x$ such that $\phi(y)= \phi(x)$ then we define
$\iota(x)=y$. Otherwise, we define $\iota(x)=x$. For every $e \in
E(G)$ such that $\phi(e) \in E(G')$, there is a unique edge $e' \in
E(G)$ such that $e' \neq e$ and $\phi(e')=\phi(e)$, and we define
$\iota(e)=e'$.  Define $\iota(e)=e$ for every $e \in E(G)$ such that
$\phi(e) \in V(G')$.

If $x \in V(G)$, $e \in E(G)$, $x \in e$  and $\phi(e) \in E(G')$
then either $\iota(x) \in \iota(e)$, or $x \in \iota(e)$. In the
second case, $m_{\phi}(x)=2$ and therefore by non-degeneracy of
$\phi$ we have $x = \iota(x)$. It follows easily from this that $\iota$ is a
morphism. Clearly $\iota \circ \iota$ is the identity map. In
particular, $\iota$ is a bijection. Therefore $\iota$ is an
involution, and $\iota$ is mixing by definition. Finally, it is easy
to see that $\phi=\iota^{\sim}$.

Now suppose $|V(G)|>2$, and let $\iota: G \rightarrow G$ be  a
mixing involution of $G$.  Denote $G/\iota$ by $G'$.  Note that
$|V(G')|\geq |V(G)|/2 >1$. Consider a vertex $x \in V(G)$, let $y =
\iota^{\sim}(x)$, and consider an edge $e'=yy'\in E(G')$. Then there
exists an edge $e=xx'$ in $G$ such that $\iota^{\sim}(e)=e' $, and
$(\iota^{\sim})^{-1}(e')=\{e, \iota(e)\}$. Therefore $|\{d \in E(G)
|x \in d,\; \iota^{\sim}(d) =e' \}|=1$ if $x \neq \iota(x)$ and
$|\{d \in E(G) |x \in d,\; \iota^{\sim}(d) =e' \}|=2$ otherwise. It
follows that $m_{\iota^{\sim}}(x)$ is well defined and positive,
and that $\sum_{\iota^{\sim}(y)=z} m_{\iota^{\sim}}(y)=2$ for every
$z \in V(G')$. Therefore, $\iota^{\sim}$ is a non-degenerate harmonic
morphism of degree two, as claimed.
\end{proof}

\medskip

The following result will be used to reduce the study of general
hyperelliptic graphs to the special case of graphs which are
2-edge-connected. 

\begin{lem}
\label{lem:ContractionEquiv}
Let $G$ be a graph, let $\overline{G}$ be the graph obtained by contracting every bridge of $G$, 
and let $\rho : G \to \overline{G}$ be the natural surjective morphism.
Then for every divisor $D \in \Div(G)$, we have $D \sim_G 0$ if and only if $\rho_*(D) \sim_{\overline{G}} 0$,
where $\rho_*(D)$ is defined as in (\ref{eq:Lowerstar}).
\end{lem}

\begin{remark}
Note that the morphism $\rho : G \to \overline{G}$ is not necessarily harmonic, c.f.~Example~\ref{ContractionExample}.
\end{remark}

\begin{proof}[Proof of Lemma~\ref{lem:ContractionEquiv}]
It suffices by induction to prove the result with $\overline{G}$ replaced by
the graph obtained by contracting a {\em single} bridge $e$.
We begin with some notation.  Let $x_1,x_2$ be the endpoints of $e$, and let $\overline{x} = \rho(x_1) = \rho(x_2)$.
Let $G_1,G_2$ be the connected components of $G-e$ containing $x_1$ and $x_2$, respectively,
and for $i=1,2$, let $\overline{G}_i = \rho(G_i)$, so that
$\overline{G} = \overline{G}_1 \cup \overline{G}_2$ and
$\overline{G}_1 \cap \overline{G}_2 = \{ \overline{x} \}$.
Note that $(x_1) \sim (x_2)$ on $G$; this follows from the observation that
$(x_1) - (x_2) = \ddiv(\chi_{G_1})$.

Let $D \in \Div(G)$. Suppose first that $D$ is a principal divisor
on $G$; we want to show that $\rho_*(D)$ is a principal divisor on
$\overline{G}$.  It suffices by linearity to consider the case where
$D = \ddiv(\chi_y)$ for some $y \in V(G)$. If $y \not\in \{ x_1,x_2
\}$, then $\rho_*(D) = \ddiv(\chi_{\rho(y)})$.  Otherwise, we have
$\rho_*(\ddiv(\chi_{x_1})) = \ddiv(\chi_{V(\overline{G}_2)})$ and
$\rho_*(\ddiv(\chi_{x_2})) = \ddiv(\chi_{V(\overline{G}_1)})$. This
proves that $\rho_*(D)$ is principal.

In the other direction, suppose that $\rho_*(D)$ is principal; we want to show that $D$ itself is principal.
By linearity, it suffices to consider the case where $\rho_*(D) = \ddiv(\chi_z)$ for some $z \in V(\overline{G})$.
If $z \neq \overline{x}$, then $\rho^{-1}(z)$ consists of a single element, and $D = \ddiv(\chi_{\rho^{-1}(z)})$.
If $z = \overline{x}$, then $\rho^{-1}(z) = \{ x_1,x_2 \}$ and using the fact that $(x_1) \sim (x_2)$ it
is easy to see that $D \sim \ddiv(\chi_{ \{ x_1,x_2 \}})$.
This proves that $D$ is principal.
\end{proof}

\begin{remark}
As alluded to in \cite[Remark~4.8]{BakerNorine},
one can use Lemma~\ref{lem:ContractionEquiv} to obtain an alternate proof of
Corollary~4.7 from \cite{BakerNorine} which does not make use of circuit theory.
\end{remark}

\begin{cor}
\label{cor:ContractionEquiv1}
Let $G$ be a graph, let $\overline{G}$ be the graph obtained by contracting every bridge of $G$, 
and let $\rho : G \to \overline{G}$ be the natural surjective morphism.
Then for every divisor $D \in \Div(G)$, we have $r_G(D) = r_{\overline{G}}(\rho_*(D))$.
\end{cor}

\begin{proof}
Let $k\geq 0$ be an integer, and let $D \in \Div(G)$.
Suppose $r(D) \geq k$, and let $\overline{D} = \rho_*(D)$.
Then for every effective divisor $E \in \Div(G)$ of degree $k$,
there exists an effective divisor $E' \in \Div(G)$ such that $D - E \sim E'$,
and thus $\overline{D} - \rho_*(D) \sim \rho_*(E')$ by
Lemma~\ref{lem:ContractionEquiv}.  Since $\rho_* : \Div(G) \to \Div(\overline{G})$
is surjective and preserves degrees and effectivity, it follows that
$r(\overline{D}) \geq k$.

Conversely, suppose $r(\rho_*(D)) \geq k$.
Then for every effective divisor $E \in \Div(G)$ of degree $k$,
there exists an effective divisor $E' \in \Div(G)$ such that $\rho_*(D) - \rho_*(E) \sim \rho_*(E')$.
By Lemma~\ref{lem:ContractionEquiv}, it follows that $D - E \sim E'$, and thus $r(D) \geq k$ as desired.
\end{proof}

\begin{cor}
\label{cor:ContractionEquiv2}
Let $G$ be a graph, and let $\overline{G}$ be the graph obtained by contracting every bridge of $G$.
Then $G$ is hyperelliptic if and only if $\overline{G}$ is hyperelliptic.
\end{cor}

\begin{proof}
This follows immediately from Corollary~\ref{cor:ContractionEquiv1} and the surjectivity of
$\rho_* : \Div(G) \to \Div(\overline{G})$.
\end{proof}

Because of Corollary~\ref{cor:ContractionEquiv2}, when studying hyperelliptic graphs there is no loss of
generality if we restrict our attention
to graphs which are $2$-edge-connected.  
And it turns out that for $2$-edge-connected graphs, there are several equivalent characterizations of what it means to be hyperelliptic:

\begin{thm}\label{t:Hyperelliptic} For a $2$-edge-connected graph $G$ of genus $g \geq 2$, the following
conditions are equivalent:
\begin{enumerate}
  \item\label{c:Condition1} $G$ is hyperelliptic.
  \item\label{c:Condition2} There exists an involution $\iota: G \rightarrow G$
  such that $G /\iota$ is a tree.
  \item\label{c:Condition3} There exists a non-degenerate
  degree two harmonic morphism $\phi$ from $G$ to a tree, or $|V(G)|=2$.
\end{enumerate}
\end{thm}
\begin{proof}
If $|V(G)|=2$ then it is easily verified that conditions (\ref{c:Condition1}),
(\ref{c:Condition2}) and (\ref{c:Condition3}) all hold.
Therefore in what follows we assume $|V(G)| > 2$.

$\mathbf{(\ref{c:Condition1}) \Rightarrow (\ref{c:Condition2})}$.  
Let $D$ be a divisor of degree $2$ on $G$ with $r(D) = 1$.
For every $x \in V(G)$, we have $|D -(x)| \neq \emptyset$ and
$\deg(D-(x))=1$. Since $G$ is $2$-edge-connected, there exists a
unique $y \in V(G)$ such that $D-(x) \sim (y)$. Define $\iota(x)=y$.

Our next goal is to define $\iota$ on $E(G)$.
Consider an edge $e=xy \in E(G)$.
If $\iota(x)=y$, we define $\iota(e)=e$. If $\iota(x) \neq
y$, then let $D_1 = (x)+(\iota(x))$ and let $D_2 = (y) + (\iota(y))$.
By the definition of $\iota$, we have $D_1 \sim D \sim D_2$. Therefore,
there exists a non-constant function $f: V(G) \rightarrow \ZZ$ such that $D_1
- D_2 = \ddiv(f)$. Let $M(f)$ be the set of all the vertices $z \in
V(G)$ for which $f(z)$ is maximal. For every vertex $z \in M(f)$, we
have

$$D_1(z) \geq (\ddiv(f))(z) = \sum_{e' = zz' \in E(G)} (f(z)-f(z')) \geq |\{e' = zz' \in E(G) \;|\; z' \in V(G)\setminus
M(f)\}|.$$

Therefore $\deg(D_1) \geq |\delta(M(f))|$,
where for $X \subseteq V(G)$ we denote by $\delta(X)$ the set of all edges of $G$
having exactly one end in $X$.
On the other hand,
$|\delta(M(f))| \geq 2$ by the $2$-edge connectivity of $G$. It follows
that $|\delta(M(f))| = 2$, and that $x, \iota(x) \in M(f)$. Analogously, we
can conclude that $f$ is minimized on $y$ and $\iota(y)$, and
therefore that $y, \iota(y) \in V(G)\setminus M(f)$. It follows that $e \in
\delta(M(f))$. Define $\iota(e)$ to be the unique edge $e^*$ such that
$\delta(M(f)) = \{e,e^*\}$.  Let $x'$ be the end of $e^*$ in $M(f)$. By
the argument above we have $D_1=(x')+(x)$. Therefore $x' =
\iota(x)$. By the symmetry between $x$ and $y$, we conclude that $e^*$
joins $\iota(x)$ and $\iota(y)$. Therefore $\iota$ is an
automorphism, and clearly $\iota \circ \iota$ is the identity.

By Lemma~\ref{l:Involution}, we know that $\phi=\iota^{\sim}$ is a
harmonic morphism. For every $x,y \in V(G /\iota)$ we have
$$\phi^*((x))= (x)+(\iota(x)) \sim D \sim (y)+(\iota(y))=
\phi^*((y)).$$ Therefore, by Theorem~\ref{thm:injthm},
we have $(x)\sim (y)$ for all $x,y \in V(G)$.
It follows from Lemma~\ref{l:allequiv} that $G /\iota$ is a tree, as
desired.

$\mathbf{(\ref{c:Condition2}) \Leftrightarrow (\ref{c:Condition3})}$.
Consider an involution $\iota$ satisfying (\ref{c:Condition2}). For
every edge $e=xy \in E(G)$ such that $x \neq \iota(y)$, the set of
edges $\{e, \iota(e)\}$ is the preimage of an edge of $G /\iota$,
and therefore forms a cut in $G$. It follows that $e \ne \iota(e)$, and therefore  $\iota$ is mixing.
The equivalence of (\ref{c:Condition2}) and
(\ref{c:Condition3}) now follows from Lemma~\ref{l:Involution}.

$\mathbf{(\ref{c:Condition3}) \Rightarrow (\ref{c:Condition1})}$.
Let $\phi: G \rightarrow T$ be a
non-degenerate harmonic morphism of degree two, where $T$ is a tree.
Let $y_0 \in V(T)$ be chosen arbitrarily and let $D :=
\phi^*((y_0))$. Then $D$ is an effective divisor of degree $2$ on $G$. 
We claim that $r(D)=1$.
Clearly, $r(D) \leq 1$. Therefore, it suffices to show that $| D -
(x)| \neq \emptyset$ for every $x \in V(G)$. Note that $(y)\sim(y')$
for every pair of vertices $y,y' \in V(T)$. Therefore $(\phi(x))
\sim (y_0)$, and by Proposition~\ref{p:PullPushFunct} we have $D \sim \phi^*((\phi(x))) \geq
m_{\phi}(x)(x)$. By since $\phi$ is non-degenerate, we have $m_{\phi}(x)>0$,
and therefore $\phi^*((\phi(x))) = (x) + (x')$ for some $x' \in V(G)$, which implies that 
$|D-(x)| \neq \emptyset$ as desired.
\end{proof}

\begin{remark}
One can use Theorem~\ref{t:Hyperelliptic} to give an alternate proof of 
Lemma~\ref{lem:hyperelledgeconn} which does not make use of 
Theorem~\ref{t:GraphAbelJacobi}.
Indeed, if $G$ is $2$-edge-connected and $|V(G)| > 2$, then 
by Theorem~\ref{t:Hyperelliptic} there is a non-degenerate harmonic morphism $\phi$ of degree $2$ from
$G$ to a tree $T$ with $|E(T)| > 0$.  If $e' \in E(T)$ and $e,\iota(e)$ are the distinct edges of $G$ mapping to $e'$
under $\phi$, then it is easy to see that $G-\{e,\iota(e) \}$ is disconnected.  Thus $G$ is not
$3$-edge-connected.
\end{remark}

It is worth stating explicitly the following fact which was established during the course of 
our proof of Theorem~\ref{t:Hyperelliptic}:

\begin{cor}
\label{cor:hyperellipticrank}
If $G$ is a $2$-edge-connected hyperelliptic graph, then
for any involution $\iota$ for which $G / \iota$ is a tree, 
we have $(x) + (\iota(x)) \sim (y) + (\iota(y))$ for all $x,y \in V(G)$.
In particular, $r((x) + (\iota (x))) = 1$ for all $x \in V(G)$.  
\end{cor}

From Corollary~\ref{cor:hyperellipticrank} and Proposition~\ref{prop:uniqueg12},
we obtain the following graph-theoretic result whose statement does not involve
harmonic morphisms at all:

\begin{cor}
\label{cor:uniqueinvolution}
If $G$ is a $2$-edge-connected graph of genus at least $2$, then there is at most one
involution $\iota$ of $G$ whose quotient is a tree.
\end{cor}

\begin{proof}
By Corollary~\ref{cor:hyperellipticrank}, if $\iota$ is such an involution then
$r((x) + (\iota (x))) = 1$ for all $x \in V(G)$.  So if $\iota$ and $\iota'$ are two
such involutions, then $(x) + (\iota(x)) \sim (x) + (\iota'(x))$ for all $x \in V(G)$
by Proposition~\ref{prop:uniqueg12}.  Thus $(\iota(x)) \sim (\iota'(x))$ for all $x \in V(G)$.
Since $G$ is $2$-edge-connected, it follows from Theorem~\ref{t:GraphAbelJacobi} that $\iota(x) = \iota'(x)$
for all $x \in V(G)$, i.e., $\iota = \iota'$.
\end{proof}

If $G$ is a $2$-edge-connected hyperelliptic graph, we call the 
unique involution $\iota$ whose quotient is a tree the {\em hyperelliptic involution} on $G$.

\begin{remark}
\label{uniquenessremark}
Corollary~\ref{cor:uniqueinvolution} is the graph-theoretic analogue of the fact that
the hyperelliptic involution on a hyperelliptic Riemann surface is unique.  
We will give another proof of Corollary~\ref{cor:uniqueinvolution} in Remark~\ref{uniquenessremark2} below.
\end{remark}

\begin{remark}
\label{involutionrankremark}
It follows from the proofs of Theorem~\ref{t:Hyperelliptic} and Corollary~\ref{cor:uniqueinvolution} that
if $G$ is a $2$-edge-connected hyperelliptic graph and $r((x)+(y))=1$ for some $x,y \in V(G)$,
then $y = \iota(x)$.
\end{remark}

As a consequence of the uniqueness of the hyperelliptic involution, we obtain the following 
corollary:

\begin{cor}
\label{cor:centralinvolution}
If $G$ is a $2$-edge-connected hyperelliptic graph with hyperelliptic involution $\iota$,
then $\iota$ belongs to the center of the group $\Aut(G)$.
\end{cor}

\begin{proof}
Let $\tau \in \Aut(G)$, and consider the automorphism $\iota' := \tau^{-1}\iota \tau$.
It is easy to check that $\iota'$ is an involution, and that $\tau$ induces an isomorphism from $G/\iota'$ to $G/\iota$, so that
$G/\iota'$ is a tree.  By Corollary~\ref{cor:uniqueinvolution}, we have $\iota' = \iota$, and therefore
$\iota$ and $\tau$ commute, as desired.
\end{proof}


\subsection{Equivalent characterizations of the hyperelliptic involution}

For a Riemann surface $X$ of genus at least $2$ and $\iota : X \to X$ an automorphism, the following are equivalent: 
(i) $X$ is hyperelliptic with hyperelliptic involution $\iota$; (ii) $\iota_* : \Jac(X) \to \Jac(X)$ is multiplication by $-1$;
(iii) $\iota^* : \Jac(X) \to \Jac(X)$ is multiplication by $-1$; 
(iv) $\iota_* : \Omega^1(X) \to \Omega^1(X)$ is multiplication by $-1$; and
(v) $\iota^* : \Omega^1(X) \to \Omega^1(X)$ is multiplication by $-1$.  
We now show that a similar characterization holds for $2$-edge-connected graphs with genus at least 2.

\begin{thm}
\label{t:hyperellautom}
Let $G$ be a $2$-edge-connected graph of genus $g \geq 2$, and let $\iota \in \Aut(G)$.
Then the following are equivalent:
\begin{enumerate}
\item $G$ is hyperelliptic with hyperelliptic involution $\iota$.
\item $\iota_* : \Jac(G) \to \Jac(G)$ is multiplication by $-1$.
\item $\iota^* : \Jac(G) \to \Jac(G)$ is multiplication by $-1$.
\item $\iota_* : \H^1(G) \to \H^1(G)$ is multiplication by $-1$.  
\item $\iota^* : \H^1(G) \to \H^1(G)$ is multiplication by $-1$.  
\end{enumerate}
\end{thm}

\begin{proof}
Since $\iota$ is a harmonic morphism of degree $1$ from $G$ to itself, $\iota_* \circ \iota^*$ is the identity map
on both $\Jac(G)$ and $\H^1(G)$.  It follows easily that $(2) \Leftrightarrow (3)$ and $(4) \Leftrightarrow (5)$.
So it suffices to prove that $(1) \Leftrightarrow (2)$ and $(1) \Leftrightarrow (5)$.

${\bf (1) \Rightarrow (2)}$.  
If $G$ is hyperelliptic with hyperelliptic involution $\iota$, then by Corollary~\ref{cor:hyperellipticrank}, 
for every $x,y \in V(G)$, we have $(x) + (\iota(x)) \sim (y) + (\iota(y))$.  
Thus $(x) - (y) \sim (\iota(y)) - (\iota(x)) = \iota_*((y) - (x))$. 
Since the group $\Div^0(G)$ is generated by divisors of the form $(x) - (y)$, it follows that $\iota_* \equiv -1$ on $\Jac(G)$.

${\bf (2) \Rightarrow (1)}$.  
If $\iota_* \equiv -1$ on $\Jac(G)$, then 
$(x) + (\iota(x)) \sim (y) + (\iota(y))$ for all $x,y \in V(G)$.
In particular, for any $x \in V(G)$, we have
$r((x) + (\iota(x))) = 1$.  Thus $G$ is hyperelliptic. 
By Remark~\ref{involutionrankremark}, $\iota$ is the hyperelliptic involution on $G$.

${\bf (1) \Rightarrow (5)}$.  
Suppose $G$ is hyperelliptic with hyperelliptic involution $\iota$, and let $\pi : G \to T$ be the corresponding quotient map 
from $G$ to a tree $T$.  If $|V(G)| = 2$, it is clear that (5) holds, so we may assume that $|V(T)| > 1$.
Let $e' \in \vec{E}(T)$ be a directed edge of $T$.
Since $\pi$ is a harmonic morphism of degree 2, 
there are two distinct directed edges $e,\iota(e)$
of $G$ mapping onto $e'$.  Let $\omega \in \H^1(G)$.
Since $T$ is a tree, we have $\H^1(T) = 0$, and therefore $(\pi_*\omega)(e') = 0$.
On the other hand, by definition we have
\[
(\pi_*\omega)(e') = \omega(e) + \omega(\iota(e))
= \omega(e) + (\iota^*\omega)(e).
\]
Since $\pi$ is surjective on oriented edges, 
it follows that $(\iota^*\omega)(e) = -\omega(e)$ for all $e \in \vec{E}(G)$ such that 
$\pi(e) \in \vec{E}(T)$.  But for $e \in \vec{E}(G)$ with $\pi(e) \in V(T)$, we have
$\iota(e) = \overline{e}$, and thus $(\iota^*\omega)(e) = -\omega(e)$ for such edges as well.
It follows that $(\iota^*\omega)(e) = -\omega(e)$ for all $e \in \vec{E}(G)$, as desired.

${\bf (5) \Rightarrow (1)}$.  
Suppose $\iota^* \equiv -1$ on $\H^1(G)$.  Then $(\iota^2)^*$ is the identity map on $\H^1(G)$, so 
$\iota$ is an involution by Proposition~\ref{prop:autinj}.
If $\iota(e) = e$ for some directed edge $e$, then
letting $\omega$ be the characteristic function of any simple cycle containing $e$, we have
\[
\omega(e) = \omega(\iota(e)) = (\iota^*\omega)(e) = -\omega(e),
\]
so that $\omega(e) = 0$, a contradiction. 
Therefore $\iota$ is mixing.
If $|V(G)| = 2$, it is easy to verify directly that (1) holds.  
So we may assume without loss of generality that $|V(G)| > 2$.
By Lemma~\ref{l:Involution}, we know that $\pi := \iota^{\sim} : G \to G' := G/\iota$ is a non-degenerate 
harmonic morphism of degree $2$.  It remains to show that $G'$ is a tree.
Since $\pi \circ \iota = \pi$, we have $\iota^*(\pi^*(\omega')) = \pi^*(\omega')$ for every $\omega' \in \H^1(G')$ by functoriality.
Since $\iota^* \equiv -1$ on $\H^1(G)$, we conclude that
$\pi^*(\omega') = -\pi^*(\omega')$, and therefore $\pi^*(\omega') = 0$, for every $\omega' \in \H^1(G')$.
But $\pi^* : \H^1(G') \to \H^1(G)$ is injective, so it follows that $\H^1(G') = 0$, i.e., $G'$ is a tree.
\end{proof}

\begin{remark}
\label{uniquenessremark2}
Combining Proposition~\ref{prop:autinj} with the proof of $(1) \Rightarrow (5)$ in
Theorem~\ref{t:hyperellautom} yields another proof of Corollary~\ref{cor:uniqueinvolution}
(i.e., of the uniqueness of the hyperelliptic involution).
\end{remark}

As an application of Theorem~\ref{t:hyperellautom}, we establish a special case of 
\cite[Conjecture~3.14]{BakerSpecialization}.
To state the result, given a graph $G$ and a positive integer $k$, we define $\sigma_k(G)$ to be the graph
obtained by replacing each edge of $G$ by a path consisting of $k$ edges.

\begin{cor}
Let $G$ be a graph, and let $k$ be a positive integer.  Then $G$ is hyperelliptic if and only if 
$\sigma_k(G)$ is hyperelliptic.
\end{cor}

\begin{proof}
By Corollary~\ref{cor:ContractionEquiv2}, we may assume without loss of generality that $G$ (and therefore $\sigma_k(G)$ as well) is a 2-edge-connected
graph of genus at least $2$.  If $G$ is hyperelliptic, then by Theorem~\ref{t:hyperellautom} there is an automorphism $\iota$ of $G$ which acts as $-1$ on
$\H^1(G)$.  Identifying $V(G)$ with a subset of $V(\sigma_k(G))$ in the obvious way induces an isomorphism between $\H^1(G)$ and $\H^1(\sigma_k(G))$, and 
it is easy to see that $\iota$ can be extended to an automorphism of $\sigma_k(G)$ which acts as $-1$ on $\H^1(\sigma_k(G))$.  Therefore
$\sigma_k(G)$ is hyperelliptic.  
Conversely, suppose that $\sigma_k(G)$ is hyperelliptic.  Then there is an automorphism $\iota'$ of $G' := \sigma_k(G)$ which acts as $-1$ on $\H^1(G')$.
By an argument similar to the proof of Proposition~\ref{prop:autinj}, it follows that $\iota'$ induces an automorphism $\iota$ of $G$ 
which acts as $-1$ on $\H^1(G)$ (the key point is that every cycle in $G'$ contains a vertex of degree at least $3$, which must belong to $V(G)$,
and which must be sent by $\iota$ to another such vertex).
Therefore $G$ is hyperelliptic as well.
\end{proof}

\medskip

\subsection{The canonical map and $3$-edge-connectivity}

We now turn to a discussion of a graph-theoretic analogue of the ``canonical map'' from a Riemann surface
to projective space.  In algebraic geometry, the following are equivalent for a Riemann surface $X$ of genus at least $2$: 
(i) $X$ is not hyperelliptic; (ii) the symmetric square $S^{(2)} : \Div_+^{2}(X) \to \Jac(X)$ of the Abel-Jacobi map 
is injective; and (iii) the canonical map $\psi_X : X \to \PP(\Omega^1(X))$ is injective.
We have already seen that the analogues of (i) and (ii) are not equivalent for $2$-edge-connected graphs 
of genus at least $2$; indeed, by Theorem~\ref{t:GraphAbelJacobi}, $S^{(2)} : \Div_{+}^{2}(G) \to \Jac(G)$ is injective if and only if
$G$ is $3$-edge-connected, and
this is a strictly weaker condition than $G$ being non-hyperelliptic (if $|V(G)| > 2$).
We now define a graph-theoretic version $\psi_G$ of the canonical map, and show that the 
analogues of conditions (ii) and (iii) for graphs are equivalent.  In other words, we will show that
$\psi_G$ is injective if and only if $G$ is $3$-edge-connected.

\medskip

Let $G$ be a $2$-edge-connected graph, and let $\H^1(G)$ be the space of harmonic $1$-forms on $G$, as defined in
\S\ref{FunctorialSection}.  We write $\PP(\H^1(G))$ for the projective space consisting of all hyperplanes 
(linear subspaces of codimension $1$) in $\H^1(G)$.  We define the {\em canonical map} 
$\psi_G : E(G) \to \PP(\H^1(G))$ by sending an edge $e \in E(G)$ to the hyperplane 
$W(e) := \{ \omega \in \H^1(G) \; : \; \omega(e) = 0 \}$.  Note that the condition $\omega(e) = 0$ is independent of
the orientation of $e$, so it makes sense to ask whether or not $\omega$ vanishes on an undirected edge.
Also, the fact that $G$ is $2$-edge-connected guarantees that $W(e) \neq \H^1(G)$, so $W(e)$ is indeed a hyperplane.

\medskip

Our main observation about the canonical map is the following proposition:

\begin{prop}
Let $G$ be a $2$-edge-connected graph.  Then the following are equivalent:
\begin{enumerate}
\item The canonical map $\psi_G : E(G) \to \PP(\H^1(G))$ is injective.
\item The map $S^{(2)} : \Div_{+}^{2}(G) \to \Jac(G)$ is injective.
\item $G$ is $3$-edge-connected.
\end{enumerate}
\end{prop}

\begin{proof}
We already know by Theorem~\ref{t:GraphAbelJacobi} that $(2) \Leftrightarrow (3)$, so it suffices to prove that 
$(1) \Leftrightarrow (3)$. 
Suppose first that $G$ is $3$-edge-connected, and let $e_1,e_2 \in E(G)$.
Since $G-\{ e_1,e_2 \}$ is connected, there is a cycle $C$ containing $e_1$ but not $e_2$.
The characteristic function $\chi_C$ of $C$ is then a flow belonging to $W(e_2)$ but not $W(e_1)$, from
which it follows that $\psi_G$ is injective.

Conversely, suppose $G$ is not $3$-edge-connected.  Then
there exist edges $e_1,e_2 \in E(G)$ such that $G-\{ e_1,e_2 \}$ is disconnected.  
It follows that any flow $\omega \in \H^1(G)$ which is non-zero on $e_1$ must also be non-zero on $e_2$.
Thus $W(e_1) = W(e_2)$, and $\psi_G$ is not injective.
\end{proof}

\begin{remark}
One can define an analogue $\psi_{G,A}$ of the canonical map for flows with values in an arbitrary abelian group $A$, and 
certain graph-theoretic assertions about $A$-flows translate nicely into statements about $\psi_{G,A}$.
For example, for $A = \ZZ/5\ZZ$,
Tutte's famous $5$-flow conjecture (c.f.~\cite[\S{X.4}, p.~348]{Bollobas})
is equivalent to the assertion that if $G$ is a $2$-edge-connected graph, 
then the image of $\psi_{G,\ZZ/5\ZZ} : E(G) \to \PP(H^1(G,\ZZ/5\ZZ))$ is
contained in an affine subspace 
(i.e., there exists a hyperplane in $\PP(H^1(G,\ZZ/5\ZZ))$ disjoint from $\psi_{G,\ZZ/5\ZZ}(E(G))$).
\end{remark}


\medskip

\subsection{Hyperelliptic graphs without Weierstrass points}

We conclude by using Theorem~\ref{t:Hyperelliptic} and the Riemann-Hurwitz formula for graphs to give
a complete characterization of all hyperelliptic graphs having no Weierstrass points.
(Graphs with no Weierstrass points are quite interesting from the point of view of
arithmetic geometry, c.f.~\cite[Corollary~4.10]{BakerSpecialization}.)

\medskip

Recall from \cite{BakerSpecialization} that, by analogy with the theory of Riemann surfaces, 
a vertex $x \in V(G)$ is called a {\em Weierstrass point} if $r(g(x)) \geq 1$.
An example is given in \cite{BakerSpecialization} of a family of graphs
of genus at least $2$ with no Weierstrass points, namely the family $B_n = B(1,1,\ldots,1)$
consisting of two vertices joined by $n\geq 3$ edges.
This is in contrast to the classical situation, in which every Riemann surface of genus at least $2$
has Weierstrass points.
(It is also proved in in \cite{BakerSpecialization} that every
{\em metric graph} of genus at least $2$ does have Weierstrass points.)

\begin{remark}
\label{InvolutionFixedPointRemark}
On a hyperelliptic Riemann surface $X$, the Weierstrass points are precisely the fixed points of the
hyperelliptic involution.  For a $2$-edge-connected graph $G$, it is easy to see that a fixed point of the
hyperelliptic involution is a Weierstrass point, and if $g(G) = 2$ then the converse also holds.  However, 
if $g(G) \geq 3$ then the converse does not always hold, as the following example shows.
Let $G$ be the hyperelliptic graph $B(3,3,3,3)$ of genus $3$, 
and let $x,y \in V(G)$ be the internal vertices of one of the edges of $G$.
Then it is not hard to verify 
that $x$ and $y$ are Weierstrass points. 
Since $\iota(x) = y$, we see that these points are not fixed by the hyperelliptic involution $\iota$ on $G$.
\end{remark}

It turns out that apart from a few
exceptions, hyperelliptic graphs almost always have Weierstrass points.
One exception is the family of graphs $B_n$ mentioned above.
Another is the family of graphs $B(l_1,l_2,l_3)$, where
$l_1,l_2,l_3$ are odd positive integers (c.f.~Example~\ref{ex:Banana}).
A third exception is the family of graphs $\Phi(l)$ described in the next paragraph.

For every integer $l \geq 1$, let the graph $\Phi(l)$ consist of
two disjoint paths $P=[x_0,x_1,\ldots,x_l]$ and $Q=[y_0,y_2,\ldots,y_l]$ of length $l$,
together with two pairs of parallel edges joining $x_0$ to $y_0$ and $x_l$ to $y_l$,
respectively.
It is easy to verify that for the unique involution $\iota : G \to G$ sending $x_i$ to $y_i$, 
the quotient graph $\Phi(l)/\iota$ is isomorphic to a path of length $l$.  
Thus $\Phi(l)$ is hyperelliptic for all $l$.

\begin{remark}
It follows from Corollary~\ref{cor:ContractionEquiv1} that
$x \in V(G)$ is a Weierstrass point if and only if $\rho(x) \in V(\overline{G})$ is
a Weierstrass point, where $\overline{G}$ is the $2$-edge-connected graph obtained by contracting every bridge of $G$.
So without loss of generality, when studying Weierstrass points on graphs it suffices to consider graphs
which are $2$-edge-connected.
\end{remark}

\begin{thm}
\label{t:HyperellipticWeierstrass}
The following are the only $2$-edge-connected hyperelliptic graphs with no Weierstrass points:

\begin{enumerate}
\item\label{c:Weierstrass2} The graph $B_n$ for some integer $n \geq 3$.
\item\label{c:Weierstrass3} The graph $B(l_1,l_2,l_3)$ for some odd integers $l_1,l_2,l_3 \geq 1$.
\item\label{c:Weierstrass4} The graph $\Phi(l)$ for some integer $l \geq 1$.
\end{enumerate}

\end{thm}

\begin{proof} Let $G$ be a $2$-edge-connected hyperelliptic graph with no Weierstrass points.
If $|V(G)| = 2$, then $G$ is isomorphic to $B_n$ for some $n \geq 3$,
so without loss of generality, we may assume that $|V(G)|>2$.
By Theorem~\ref{t:Hyperelliptic}, there exists a non-degenerate degree $2$
harmonic morphism $\phi: G \rightarrow T$ for some tree $T$ with $|V(T)| > 1$.
Note that for every $t \in V(T)$ we have $r(\phi^*((t))) = 1$. If $m_{\phi}(x)=2$ for some
$x \in V(G)$, then $x$ is a Weierstrass point, as $r(g(x)) \geq r(2(x)) = r(\phi^*(\phi(x)))=1$.
Therefore, we may assume without loss of generality that $m_{\phi}(x)=1$ for every $x \in V(G)$, so
that every $t \in V(T)$ has exactly two preimages under $\phi$.
By the Riemann-Hurwitz formula for graphs (Theorem~\ref{t:RiemannHurwitz}), 
we have $\sum_{x \in V(G)}v_{\phi}(x) = 2g+2$.

Consider a vertex $t \in V(T)$ with $\deg(t)=1$. Let $\phi^{-1}(t)=\{x,x'\}$, and let $x''$ be the unique neighbor
of $x$ in $V(G) \setminus \phi^{-1}(t)$ (which is well-defined since $m_{\phi}(x) = 1$). 
It is easy to see that there are $v_{\phi}(x) + 1$ edges incident to $x$, namely 
the $v_{\phi}(x)$ vertical edges connecting $x$ to $x'$ and the horizontal edge connecting $x$ to $x''$.
Also, since $G$ is $2$-edge-connected, we have $\deg(x) \geq 2$, so by (\ref{e:DegreeIdentity}) we know that $v_{\phi}(x) \geq 1$.
It follows that $$(v_{\phi}(x)+2)(x) \sim (x) + v_{\phi}(x)(x') + (x'') \geq (x)+(x') =  \phi^*(t).$$
Therefore $r((v_{\phi}(x)+2)(x)) \geq 1$, and $x$ is a Weierstrass point of $G$ if $v_{\phi}(x) \leq g - 2$.
Thus we may assume without loss of generality that $v_{\phi}(x) \geq g - 1$ for every $x \in V(G)$ such that $\deg(\phi(x))=1$. Let
$k := |\{t \in V(T) \; | \; \deg(t)=1\}| \geq 2$. We have $2g+2 = \sum_{x \in V(G)}v_{\phi}(x) \geq 2k(g-1) \geq 4(g-1).$
It follows that either $g =2$ and $k \leq 3$, or else $g = 3$ and $k=2$.  In the latter case, 
$T$ is a path and $v_{\phi}(x)=0$ for every $x \in V(G)$ such that  $\deg(\phi(x))>1$.

If $g=2$, it is easy to see that $G$ must be isomorphic to the graph $B(l_1,l_2,l_3)$ for some integers
$l_1,l_2,l_3 \geq 1$, and if $l_i$ is even for some $i \in \{1,2,3\}$ then the middle vertex of the path of length $l_i$
is a Weierstrass point by Example~\ref{ex:Banana}. If $g = 3$, then by the above we have $v_{\phi}(x)=2$ 
for every $x \in V(G)$ such that $\deg(\phi(x))=1$. 
It follows easily that $G$ is isomorphic to the graph $\Phi(|E(T)|)$.

\medskip

It remains to show that if $G$ is one of the graphs in (\ref{c:Weierstrass2}),
(\ref{c:Weierstrass3}) or (\ref{c:Weierstrass4}), then it has no
Weierstrass points. 
We start by considering the case when $G$ satisfies (\ref{c:Weierstrass2}), 
i.e., $|V(G)|=2$ and $|E(G)|=n \geq 3$.
Let $V(G)=\{x,y\}$. By Theorem~\ref{t:OrderDivisor} we have
$r((n-1)(x) - (y)) = -1$. Therefore $r(g(x))=r((n-1)(x)) \leq r((n-1)(x) - (y)) + 1 =0$.
It follows that $x$ is not a Weierstrass point, and by symmetry neither is $y$.

Now suppose that $G$ is of the form (\ref{c:Weierstrass3}), so that $g=2$.
As each $l_i$ is odd, the hyperelliptic involution $\iota$ on $G$ has no fixed points
(since, in the notation of Example~\ref{ex:Banana}, we have $\iota(x) = y$ and $\iota(z_i)=z_{l-i}$).
Therefore $G$ has no Weierstrass points by Remark~\ref{InvolutionFixedPointRemark}.

Finally, suppose that $G$ is of the form (\ref{c:Weierstrass4}), i.e.,
that $G$ is isomorphic to the graph $\Phi(l)$ for some integer $l
\geq 1$. Let the vertices of $G$ be labeled as in the definition of
$\Phi(l)$. By symmetry, it suffices to prove that $r(3(x_i))=0$ for
every integer $i$ such that $0 \leq i \leq l$. Suppose first that $l
\leq 3i \leq 2l$. Then $3(x_i) \sim (x_0)+(x_l)+(x_{3i-l})$.
Consider the following linear order $<$ on $V(G)$:
\[
y_0 < \cdots < y_l < x_0 < \cdots < x_{3i-l-1} < x_l < \cdots < x_{3i-l+1} < x_{3i-l}.
\]
The divisor associated to this order is equal to $(x_0)+(x_l)+(x_{3i-l})
- (y_0)$. It follows that $r(3(x_i))=r((x_0)+(x_l)+(x_{3i-l}))=0$ in
this case. Suppose now that $3i < l$. Then we have
$$ 3(x_i) \sim (x_{3i})+2(x_0) \sim (x_{3i+1})+2(y_0) \sim (x_l) +(y_0) + (y_{l-3i-1}).$$
This time, we consider the following linear order $<$ on $V(G)$:
\[
y_l < x_l < \cdots < x_0 < y_0 < \cdots < y_{l-3i-2} < y_{l-1} < \cdots < y_{l-3i} < y_{l-3i-1}.
\]
The divisor associated to this order is equal to $(x_l)
+(y_0) + (y_{l-3i-1}) - (y_l)$. Again it follows that $r(3(x_i))=0$.
The last remaining case, where $3i >2l$, follows by symmetry from the case $3i < l$.
\end{proof}

\begin{remark}
It would be interesting to characterize {\em all}
2-edge-connected graphs $G$ having no Weierstrass points.  We do not at present know of any examples in which $G$ is not hyperelliptic.
\end{remark}

\bibliographystyle{alpha}
\bibliography{hyperelliptic}
\end{document}